\documentclass[a4paper,12pt]{amsart}

\usepackage[all]{xy}
\usepackage{amssymb}

\usepackage{t1enc}

\def\Bbb{\mathbb}

\def\cD{\mathcal{D}}
\def\cR{\mathcal{R}}
\def\cV{\mathcal{V}}
\def\cB{\mathcal{B}}
\def\cN{\mathcal{N}}
\def\cU{\mathcal{U}}
\def\cL{\mathcal{L}}
\def\N{\mathbb{N}}
\def\R{\mathbb{R}}

\def\C{\mathbb{C}}
\def\F{\mathbb{F}}

\def\be{\beta}
\def\al{\alpha}
\def\la{\lambda}
\def\ka{\kappa}

\def\De{\Delta}

\def\bx{\mbox{\boldmath $x$}}
\def\bD{\mbox{\boldmath $\cD$}}

\newcommand{\id}{\operatorname{id}}
\newcommand{\Spec}{\operatorname{Spec}}
\newcommand{\End}{\operatorname{End}}
\newcommand{\Min}{\operatorname{Min}}
\newcommand{\Max}{\operatorname{Max}}

\newtheorem{theorem}{Theorem}[section]
\newtheorem{lemma}[theorem]{Lemma}
\newtheorem{proposition}[theorem]{Proposition}
\newtheorem{corollary}[theorem]{Corollary}

\theoremstyle{remark}
\newtheorem{remark}[theorem]{\rm\bf Remark}
\newtheorem*{definition*}{\rm\bf Definition}

\newcommand{\nn}[1]{(\ref{#1})}

\def\endrk{\hbox{$|\!\!|\!\!|\!\!|\!\!|\!\!|\!\!|$}}

\def\sideremark#1{\ifvmode\leavevmode\fi\vadjust{\vbox to0pt{\vss
 \hbox to 0pt{\hskip\hsize\hskip1em
 \vbox{\hsize3cm\tiny\raggedright\pretolerance10000
  \noindent #1\hfill}\hss}\vbox to8pt{\vfil}\vss}}}

\def\idx#1{{\em #1\/}}

\author{A. Rod Gover and Josef \v Silhan}
\email{gover@math.auckland.ac.nz} \title{Commuting linear operators
  and decompositions; applications to Einstein manifolds}
\begin{document}

\address{ARG: Department of Mathematics\\
  The University of Auckland\\
  Private Bag 92019\\
  Auckland 1\\
  New Zealand} \email{gover@math.auckland.ac.nz}
\address{JS: Eduard \v{C}ech center \\ 
Department of Algebra and geometry \\
Masaryk University \\
Jan\'a\v{c}kovo n\'am. 2a \\
602 00, Brno\\
Czech Republic} \email{silhan@math.muni.cz}

\maketitle

\pagestyle{myheadings}
\markboth{Gover \& \v Silhan}{Decomposition Theorems and Einstein manifolds}

\begin{abstract}
For linear operators which factor $P=P_0P_1 \cdots P_\ell$, with
suitable assumptions concerning commutativity of the factors, we
introduce several notions of a decomposition. When any of these hold
then questions of null space and range are subordinated to the same
questions for the factors, or certain compositions thereof. When the
operators $P_i$ are polynomial in other commuting operators then we
show that, in a suitable sense, generically factorisations algebraically yield
decompositions. In the case of operators on
a space over an algebraically closed field this boils down to
elementary algebraic geometry arising from the polynomial formula for
$P$. Applied to operators $P$ polynomial in single  other operator $\cD$
this shows that the solution space for $P$ decomposes directly into a
sum of generalised eigenspaces for $\cD$.  We give universal formulae
for the projectors administering the decomposition. In the generic
setting the inhomogenous problems for $P$ reduce to an equivalent
inhomogeneous problem for an operator linear in $\cD$. These results
are independent of the operator $\cD$, and so provide a route to
progressing such questions when functional calculus is
unavailable. Related generalising results are obtained as well as a
treatment for operators on vector spaces over arbitrary fields.  We
introduce and discuss symmetry algebras for such operators.  As a
motivating example application we treat, on Einstein manifolds, the
conformal Laplacian operators of Graham-Jenne-Mason-Sparling.
\end{abstract}

\section{Introduction}
\newcommand{\bF}{\mathbb{F}} \newcommand{\bC}{\mathbb{C}}
\newcommand{\Proj}{\operatorname{Proj}}
\renewcommand{\Pr}{\operatorname{Pr}}

A motivating algebraic question is as follows.  For $\cV$ a vector
space, $\cD:\cV\to \cV$ an arbitrary linear operator, and $P:\cV\to
\cV$ a linear operator which is polynomial in $\cD$, then what do we
know about the solution space for $P$ in terms of the generalised
eigenspaces of $\cD$?  The question is obviously most interesting when
$\cV$ is infinite dimensional. In fact we want to treat this, and
related questions, uniformly without using any information about the
operator $\cD$ or the vector space $\cV$.  Obviously any gains in this
direction are particularly important in settings where functional
calculus is unavailable, but they also provide a potentially important
first simplifying step even when there is access to functional calculus. 

In the case that the field involved is algebraically closed we obtain a
complete answer to the question above.
\begin{theorem} \label{fundthm}
Let $\cV$ be a vector over an algebraically closed field $\bF$. 
Suppose that $\cD$ is a linear endomorphism on
$\cV$, and $P=P[\cD]:\cV\to \cV$ is a linear operator polynomial in
$\cD$.  Then the solution space $\cV_P$, for $P$, admits a canonical
and unique direct sum decomposition
\begin{equation}
\cV_P=\oplus_{i=0}^\ell \cV_{\lambda_i}~,
\end{equation}
where, for each $i$ in the sum, $\cV_{\lambda_i}$ is the solution
space for $(\cD+\lambda_i)^{p_i}$ ($p_i\in \mathbb{Z}_{\geq 0}$) with 
$-\lambda_i\in \bF$ a multiplicity $p_i$  solution of the
polynomial equation $P[x]=0$. The projection $\Proj_i: \cV_P\to
\cV_{\lambda_i}$ is given by the universal formula \nn{projform:gen}.
\end{theorem}

\noindent The cross reference
\nn{projform:gen} refers an explicit formula given in the next
section. If $u\in \cV$  
satisfies 
\begin{equation}\label{power}
(\cD+\lambda )^p u=0
\end{equation} 
and is non-zero then we shall term $u$ a {\em generalised eigenvector}
for $\cD$ corresponding to the {\em generalised eigenvalue}
$-\lambda$. Using this language a partial paraphrasing of Theorem
\ref{fundthm} is that the solution space for $P$ is a direct sum of
generalised eigenspaces for $\cD$.  The Theorem above is an immediate
corollary of Theorem \ref{pf0:gen}; for the case that $P$ is given as
a fully factored expression, this states the situation for $\cV$ over
an arbitrary field. Related eigenspace/eigenspectral results
follow, see Corollary \ref{spectralthm}.

We may also consider inhomogeneous problems $P u=f$. In the case that $\bF$ is
algebraically closed then, by rescaling, this boils down to a problem 
for an operator of the form 
\begin{equation}\label{Popform}
P u := (\cD+\lambda_0)^{p_0}(\cD+\lambda_1)^{p_1} \cdots (\cD+\lambda_\ell)^{p_\ell}u .
\end{equation}
\begin{theorem}\label{inh-first}
  Let $\cV$ be a vector space over a field $\bF$.  Suppose that
  $\cD:\cV\to \cV$ is a non-trivial linear endomorphism and consider
  $P:\cV\to \cV$ as in \nn{Popform} with the $\lambda_i\in \bF$
  mutually distinct and for $i=0,\cdots ,\ell$, $p_i\in {\Bbb Z}_{\geq
    1}$.  Let us fix $f\in \cV$.  There is a 1-1 relationship between
  solutions $u\in \cV$ of $P u=f$ and solutions $(u_0,\cdots ,u_\ell)\in \oplus^{\ell+1}\cV$
  of the problem
\begin{equation} 
(\cD+\lambda_0)^{p_0}u_0=f, \cdots , (\cD+\lambda_\ell)^{p_\ell}u_\ell =f. 
\end{equation}
\end{theorem}
\noindent In Theorem \ref{inh} we give the explicit
transformation between the two problems. 

Given a linear operator $\cD':\cV\to \cV$, let us write $\cR(\cD)$ to
denote the image of $\cD'$ and $\cN(\cD)$ the kernel of $\cD'$.  We may
summarise part of the key information in Theorem \ref{inh-first} and
Theorem \ref{fundthm} (or more accurately Theorem \ref{pf0:gen}) by
 the following.
\begin{corollary} For $P:\cV\to \cV$ as in \nn{Popform} we have
$$
\cR(P)=\bigcap^{i=\ell}_{i=0}\cR((\cD+\la_i)^{p_i}), \quad
\cN(P)=\bigoplus^{i=\ell}_{i=0}\cN((\cD+\la_i)^{p_i}).
$$
\end{corollary}

Taking the special case of $D$ being $\frac{d}{dx}$ acting on the
smooth functions of $\mathbb{R}$ the above recovers much of the
standard theory of constant coefficient linear ordinary differential
equations. Evidently these aspects generalise to arbitrary operators
$\cD$.  In fact the above results are just examples from a very
general context (not a priori related to polynomials) in which we develop
considerable theory as below.

Section \ref{gensec} introduces various notions of a decomposition for
linear operators $P$ that factor $P=P_0P_1\cdots P_\ell$, and where
the factors mutually commute. Each decomposition is based on some
level of invertibility; not invertibility of $P$, nor the components
individually but rather of the system $(P_0, P_1,\cdots , P_\ell)$.
This is given initially in terms of identities satisfied by ``relative
inverses'', see \nn{iddec} and \nn{al-iddec}. However a more intuitive
picture may be obtained by diagrams as follows. One may construct a
canonical complex from the operators $P_0, P_1, \cdots ,P_\ell$ (see
the diagrams \nn{3comp} and \nn{res}). This is the Koszul (cochain)
complex for the system $(P_0,\cdots ,P_\ell)$. In each case the
required invertibility means that a certain class of subcomplexes of
this is exact. The latter is described in Section \ref{compare}. More
than this see Theorem \ref{Qfree}. This shows that, remarkably, making
only assumptions concerning the exactness of certain subcomplexes of
the full Koszul complex we recover almost the entire
information of a class of decompositions.

The basic Koszul complex which underlies several of our constructions
is also a central tool in the treatment of certain spectral systems
for commuting operators: the so-called Taylor spectrum \cite{T,Esc},
and the related split spectrum. See \cite{M}, and references therein,
for further discussion. Operators polynomial in another operator, as
above, have also been subject of spectral theory, e.g.\
\cite{MCrelle}.  While we believe there is considerable scope to
develop spectral theory based around our discussion, this will be
deferred to a later treatment. Our current focus is the use of purely
algebraic considerations which may be applied rather universally. In
particular at no point do we need a Banach structure on the vector
spaces or (spaces of) operators involved. We include some minor
comments concerning eigenspectrum and eigenspaces.

For any $P$ admitting a decomposition, of some type, and for any
inhomogeneous problem, we are able to describe completely and
explicitly the structure of the solution space in terms of data for
the component operators or certain products thereof. In particular the
results above generalise immediately, see Theorem \ref{inhg},
Corollary \ref{maincor}, Theorem \ref{winhg} and Corollary
\ref{wmaincor}.  It is meaningful to say that the required
invertibility for the system $(P_0, P_1,\cdots , P_\ell)$, in order to
obtain some decomposition, is very weak (cf.\ Chapter IV, Theorem 4
\cite{M}). 
In fact for linear operators polynomial in commuting
endomorphisms $\cD_0,\cD_1,\cdots ,\cD_k$, via elementary algebraic
geometry we show that it is attained generically.  In any specific
case, over an algebraically closed field, establishing any of the
class of decompositions boils down to verifying that collections of
algebraic varieties determined by combinations of the factors have no
common point.  See Theorem \ref{geomst}. For example constant
coefficient inhomogeneous linear partial differential equations may
generically be reduced to equivalent lower order equations using
Theorem \ref{winhg}, and specific problems are practically treatable.
We should also point out that for operators polynomial in commuting
operators the decompositions we study and obtain are constructed by
purely algebraic means. This means the results we obtain are
universal; they are independent of the operators $\cD_0,\cD_1,\cdots
,\cD_k$. The ``relative inverses'' are given by polynomials in the
same operators $\cD_0,\cD_1,\cdots ,\cD_k$. (So for example if the
$\cD_i$ are differential operators then the entire theory is within
the category of differential operators polynomial in these.) Regarding
the weakness of the relative invertibility conditions see also Chapter
IV, Theorem 4 of \cite{M}. 

For operators $P$ polynomial in a single other operator $\cD$, as
above, a polynomial factorisation of $P$ is generically the strongest
(non-trivial) form of decomposition. Thus, and for other reasons, the
strongest results are obtained in this setting. Some are summarised
above and Section \ref{mainsec} develops the full theory.

A non-trivial application for some of these ideas is the study of
differential operators polynomial in the Laplacian.  Problems of this
nature arise in differential geometry and, in particular, in the study
of conformal Laplacian type operators.  One of the simplest examples
is the conformal Laplacian $Y$. This curvature modification of the
usual Laplacian $\Delta$ is, in a suitable sense, conformally
invariant and its importance was observed early last century, see e.g.
\cite{dirac}. Paneitz constructed a conformal operator with principal
part $ \Delta^2$ \cite{Pan} and then cubic operators are due to T.\
Branson and V.\ W\"unsch.  Later Graham, Jenne, Mason and Sparling
(GJMS) \cite{GJMS} extended these results to a very large family which
in odd dimensions, for example, includes conformal Laplacian operators
of all even orders. Recently this family has been seen to play a deep
role in Riemannian, pseudo-Riemannian and conformal geometry. For
example the operators have a central role in the geometry of the
asymptotically hyperbolic Einstein-Poincar\'e metric which underlies
the AdS/CFT correspondence of physics, see e.g.\ \cite{FeffGr01,GrZ2}.
In another direction the GJMS operators control the equations for the
prescription of Branson's $Q$-curvature, and also the prescription of
the non-critical $Q$-curvatures \cite{tomsharp,DM}. These problems
generalise the celebrated Yamabe problem (see \cite{schoen} and
references therein) of scalar curvature prescription. It was shown in
\cite{GoEinst} that on Einstein manifolds the $Q$ curvature and the
non-critical $Q$-curvature are constant.  In both cases this result is
related to another result in \cite{GoEinst}, namely that on Einstein
manifolds the GJMS operators are given by factored polynomials in the
Laplacian. In section \ref{eins} we will use the Theorems above with
this polynomial factorisation to discuss in any signature, and on any
non-Ricci-flat Einstein manifold, the relationship of the solution space the
GJMS operators to the spectrum and solution space of the conformal
Laplacian operator. Via Theorem \ref{inh-first} the differential order
$2k$ inhomogeneous problem $P_{k}u=f$ for these operators may be
reduced explicitly to an equivalent second order (Laplacian) problem
of the form
$$
(\Delta+ \mbox{\boldmath$\lambda$})\mbox{\boldmath$u$}=f
$$ where $\mbox{\boldmath$\lambda$}:={\rm diag}(\la_1,\cdots,\la_k)$
(with the $\la_i$ given explicitly in terms of the scalar curvature),
$\mbox{\boldmath$u$}= {\rm Transpose}(u_1,\cdots ,u_k)$ and $f$ here
means ${\rm Transpose}(f,\cdots ,f)$, see Proposition \ref{inhPk}. (In
fact, as commented after Proposition \ref{inhPk}, by enlarging the
space on which operators may act, the problems in most cases may 
be reduced in the same spirit to equivalent first order problems.)
This may have applications in the understanding of $Q$-prescription on
conformally Einstein manifolds as such inhomogeneous problems may be
viewed as model linear problems for the true (non-linear) prescription
problems.  In fact the prescription problems involve equations of the
form $P_k u= h(u) f$, for $h$ a suitable function of $u$ (in many
cases simply $h(u)=\mbox{constant}.u^p $ for a suitable power $p$).
The tools of Theorem \ref{inh} still apply when we have a
non-linearity of this type and so such equations reduce to $ (\Delta+
\mbox{\boldmath$\lambda$})\mbox{\boldmath$u$}= h(
\mbox{\boldmath$b$}\mbox{\boldmath$u$})f, $ on non-Ricci-flat Einstein
manifolds, where $\mbox{\boldmath$b$}$ is a row matrix of constants
determined by the scalar curvature.  In yet another direction the
tools of section \ref{mainsec} show that on Einstein manifolds, of any
signature, the eigenspaces and eigenvalues of the $P_{k}$ all arise
from generalised eigenspaces and generalised eigenvalues of $Y$. This
is just a special case of the general result in Corollary
\ref{spectralthm}, and should have application in the representation
theory of the orthogonal groups $SO(p+1,q+1)$ (which, modulo issues of
covering, act as the conformal group on appropriate products of
spheres). The case of conformal Laplacians as discussed here is just
an example application. There are many other settings where these
ideas apply. For example parallel to the theory of conformal
Laplacians there is a theory of sub-Laplacians in CR geometry
\cite{GoGr}. Product manifolds yield commuting operators (such as
Laplacians of the components) and so the machinery of Section
\ref{multivar} is ready for these.
 
Recently there has been a growth in interest in the very old problem
concerning the so-called symmetries and symmetry algebras of Laplacian
type operators, see \cite{MikEsym,EL} and references therein.  Such
symmetry operators play a central role in separation of variables
techniques for the solution of the Laplacian operators involved. In
section \ref{sym} we introduce symmetry algebras which generalise this
notion to a large class of linear operators. 
Using this we obtain, for example,
 general results relating the symmetry
algebra for a linear endomorphism $\cD$ and that of a second operator
$P$ polynomial in $\cD$. See in particular Theorem \ref{genweak}. In
section \ref{natweak} Einstein manifolds are once again used to
illustrate these ideas in a concrete setting.

Finally we point out that the theory of polynomial operators may also
be applied to large classes of differential operators which are not
simply polynomial in another operator $\cD$. This is achieved by, for
example, composing operators which do have the form $P[\cD]$ (i.e.\
$P$ polynomial in a differential operator $\cD$) with other suitable
 differential operators. For example large classes of conformally
 invariant operators on tensor and spinor fields arise this way
 \cite{BrGodeRham,JoSth}. This will be taken up elsewhere.

The authors are grateful to John Butcher, Andreas \v Cap, Mike
Eastwood, V. Mathai, Paul-Andi Nagy and Jan Slov\'ak for helpful
discussions.  The first author would like to thank the Royal Society
of New Zealand for support via Marsden Grant no. 06-UOA-029.
The second author was supported from the Basic Research Center no.\ LC505 
(Eduard \v{C}ech Center for Algebra and Geometry) of Ministry of Education, 
Youth and Sport of Czech Republic.

\section{The general setup} \label{gensec}

Our study here will concern linear operators given by a composition
$P=P_0P_1\cdots P_\ell$ where the factors mutually commute.  In the
case that each factor is invertible then the essential properties of
$P$ are given simply in terms of the factors. Otherwise the situation
is significantly more complicated. Here we explore part of the latter
domain.  In general there are relative qualities of linear operators
$P_0$ and $P_1$ that dramatically affect the nature of the composition
$P_0P_1$.  As a very simple example, and a case in point, one may
compare $(\frac{d}{dx}+\la)(\frac{d}{dx}+\mu)$, with $\la \neq \mu$,
to $(\frac{d}{dx}+\la)^2$. As operators on the line these are rather
different beasts.  These observations in part motivate considering the
following class of linear operators.

\subsection{Decompositions of linear operators} \label{gendecom}

Let $\cV$ denote a vector space over a field $\bF$. Suppose
that $P:\cV\to \cV $ is a linear operator which may be expressed as a
composition 
\begin{equation}\label{compP}
P=P_0P_1\cdots P_\ell
\end{equation} where the linear operators $ P_i:\cV\to \cV$,  $ i=0,\cdots
,\ell $, have the following properties:\\ there exist linear
operators $Q_i:\cV\to \cV$, $i=0,1,\cdots, \ell$, that give a
decomposition of the identity,
\begin{equation}\label{iddec}
id_V=Q_0P^0+\cdots+Q_\ell P^\ell~,
\end{equation} 
where $P^i:=\Pi_{i\neq j=0}^{j=\ell} P_i, i=0,\cdots ,\ell$;
and the $P_i$s and the $Q_j$s are mutually commuting to the extent
\begin{equation}\label{commut}
P_iP_j=P_jP_i, \quad\mbox{and}\quad P_iQ_j=Q_jP_i,\quad i,j\in\{0,\cdots,\ell\}~;
\end{equation}  
When $\ell=0$ this may be viewed to hold trivially. 
For other cases we shall use the
following terminology.  

\begin{definition*} For a linear operator $P:\cV\to \cV$ an
expression of the form \nn{compP} will be said to be a 
{\em decomposition} 
of $P$ if the factors $P_i$, $i=0,\cdots ,\ell$,
satisfy the conditions just described and $\ell\geq 1$.  
\end{definition*}

Note that if one of the factors $P^i$ is invertible (from both sides)
then we have \nn{iddec} immediately.  In general requiring the
identity \nn{iddec} is a significantly weaker requirement. It states for example
that the operator $(P^0,\cdots,P^\ell):\oplus_{i=0}^\ell\cV\to\cV $
has a right inverse given by the operator $(Q_0,\cdots ,Q_\ell):\cV\to
\oplus_{i=0}^\ell\cV$. (This statement also holds if we swap the roles
of the $P_i$s and the $Q_i$s.) 

\vspace{1ex}

We next observe that the identity \nn{iddec} controls a decomposition
of the null space $\cN(P)$.
\begin{lemma}\label{projlem}
For each $i\in\{0,\cdots ,\ell \}$, we have 
$$
Q_iP^i:\cN(P)\to \cN(P_i)
$$ 
and this is a projection.
\end{lemma}
\begin{proof} Since $P_iQ_iP^i=Q_iP$ it is clear that 
$Q_iP^i:\cN(P)\to \cN(P_i)$. Then on $\cV$, and hence in particular on
$\cN(P)$, we have the identity given by \nn{iddec}. But for $j\neq i$, $P_i$ is a factor of $P^j$ and hence 
$P^j $ annihilates $\cN(P_i)$. So $Q_iP^i$ is the identity on $\cN(P_i)$.
\end{proof}
\noindent For convenience we will often use $\cV_P$ to denote
the null space of a linear  operator $P$ on $\cV$, so e.g.\ we may write 
$Q_iP^i:\cV_P\to \cV_{P_i}$.

We consider now the inhomogeneous problem $P u=f $. Of course the
solution space is the affine subspace in $\cV$ obtained by translating
$\cV_P$ (the solution space for the linear problem) by any single
``particular'' solution to $P u =f$. It turns out that, by applying
\nn{iddec} to $\cap_0^\ell \cR(P_i)$, we can decompose
the inhomogeneous problem to a simpler inhomogeneous problem in a way
that generalises the treatment of the homogeneous cases.
\begin{theorem}\label{inhg}
Let $\cV$ be a vector space over a field $\bF$  
and consider $P:\cV\to \cV$ as in \nn{compP} with 
the factorisation there giving  a decomposition, i.e.\ \nn{iddec} and 
\nn{commut}  hold.
Let us fix $f\in \cV$.  There is a 1-1 relationship
between solutions $u\in \cV$ of $P u=f$ and solutions $(u_0,\cdots
,u_\ell)\in\oplus^{\ell+1}\cV$ of the problem
\begin{equation}\label{multiprobg} 
P_0 u_0=f, \cdots , P_\ell u_\ell =f. 
\end{equation}

Writing $\cV_P^f$ for the solution space of $P u=f$ and (for
$i=0,\cdots ,\ell$) $\cV^f_{i}$ for the solution space of 
$P_i \tilde{u}=f$.
The map $F:\cV^f_P\to \times_{i=0}^\ell\cV^f_{i} $ is given by 
$$ u\mapsto (P^0 u, \cdots ,P^\ell u ) ~,
$$
with inverse $B:\times_{i=0}^\ell\cV^f_{i}\to \cV_P^f$ given by
$$
(u_0,\cdots, u_\ell)\mapsto \sum_{i=0}^{i=\ell} 
Q_i u_i~.
$$
On $\cV$ we have $B\circ F=id_{\cV}$, while on the affine space 
$\times_{i=0}^\ell\cV^f_{i}$ we have
$F\circ B =\id_{\times_{i=0}^\ell\cV^f_{i}}$.
\end{theorem}
\begin{proof} Suppose $Pu =f$. Then $P_iP^i u=Pu=f$ and so $Fu$ is a solution of \nn{multiprobg}. For the converse suppose that $(u_0,\cdots ,u_\ell)$ is a solution of \nn{multiprobg} and write $u:=  \sum_{i=0}^{i=\ell} 
Q_i u_i$. 
Then 
$$
\begin{aligned}
Pu= \sum_{i=0}^{i=\ell}PQ_i u_i =\sum_{i=0}^{i=\ell}Q_i P^i P_i u_i=\sum_{i=0}^{i=\ell}Q_i P^i f =f
\end{aligned}
$$
where finally we have used \nn{iddec}.

It remains to establish the final assertion. By construction $B \circ
F=\sum_{i=0}^{i=\ell} Q_i P^i $ and so $B\circ F=id_{\cV}$ is just the
identity \nn{iddec}. (Then in particular $B\circ F=id_{\cV^f_P}$.)
Next we calculate $F\circ B$ on
$\times_{i=0}^\ell\cV^f_{i}$. For the $k^{\rm th}$-component, we have 
$[FB(u_0,\cdots ,u_\ell)]_k$ given by 
$$
P^k  \sum_{i=0}^{i=\ell} 
Q_iu_i~.
$$ Using the commutativity of terms, and that
$P_iu_i=f$, this gives
$$\begin{aligned} 
\Big( \sum_{k\neq i=0}^{i=\ell}
Q_i  \prod_{i,k\neq m
 =0}^{m=\ell} P_m  f \Big)
+ Q_k P^k u_k
\end{aligned}
$$ 
Now using \nn{iddec} and then
$P_ku_k=f$, we obtain for the last term, 
$$\begin{aligned} 
Q_k P^k u_k = u_k - 
 \sum_{k\neq i=0}^{i=\ell} Q_i 
 \prod_{i,k\neq m
 =0}^{m=\ell} P_m  f ~.
\end{aligned}
$$
Thus 
$$
[FB(u_0,\cdots ,u_\ell)]_k=u_k~,
$$
for any $k\in \{0,\cdots ,\ell\}$ and we conclude that $FB$ is the identity on 
$\times_{k=0}^\ell \cV_{k}^f$.
\end{proof}

For any operator of the form \nn{compP}, with the $P_i$ mutually
commuting, we obviously have $+_{i=0}^{i=\ell}\cN(P_i)\subseteq \cN(P)$
and $\cR(P)\subseteq \cap _{i=0}^{i=\ell} \cR(P_i)$. From the above we
see that if \nn{iddec} holds then these containments are equalities.
In summary we have the following.
\begin{corollary}\label{maincor}
For $P:\cV\to \cV$, with \nn{compP} giving a decomposition, we have
$$
\cR(P)=\bigcap^{i=\ell}_{i=0}\cR(P_i), \quad
\cN(P)=\bigoplus^{i=\ell}_{i=0}\cN(P_i).
$$
The decomposition of $\cN(P)$ is given by the identity
$$
id_{\cN(P)}=\sum_{i=0}^{i=\ell} \Proj_i
$$ where, for each $i\in\{0,\cdots ,\ell\}$, $\Proj_i:\cN(P)\to
\cN(P_i)$ is the projection given by the restriction of $Q_iP^i$ from
\nn{iddec}.
\end{corollary}

In section \ref{mainsec} we shall show that operators polynomial in a
single other operator generically admit decompositions that may
obtained algebraically and explicitly.  The explicit formulae for the
$Q$'s (in the identity \nn{iddec}) are given in terms of the basic
data of the factorisation $P_0P_1\cdots P_\ell$.  Applications treated
in Sections \ref{sym} and \ref{eins} then show that the decompositions
are a powerful tool. Before we take these directions we study the
algebraic structures underlying a decomposition and this leads to
results which show that the decompositions are a special case of a
rather general theory with a vastly broader scope for development and
applications.

\subsection{Relative invertibility, and operator resolutions} \label{complexes}
We first  shed some light on the above constructions.  

In relation to the identity  \nn{iddec}, suppose that we have
linear endomorphisms $P_0,P_1$, on a vector space $\cV$, and there exist 
further endomorphisms $Q_1,Q_0$ so that
$$
id_\cV= Q_1P_0+Q_0P_1.
$$ Then clearly $P_0$ is injective on the null space of $P_1$ and is
an invertible endomorphism of $\cN(P_1)$ space if, for example, the
commutativity relations \nn{commut} hold.
So the displayed identity manifests what
we might call {\em relative invertibility} of the operators $P_0$ and
$P_1$. Evidently 
we may solve $P_0 u=f$ for $f\in \cN(P_1)$. This is a consequence of the 
fact that the short complex associated to the system, 
\begin{equation}\label{seq1}
0\to \cV \stackrel{\left(\begin{array}{c} \mbox{\scriptsize{$P_0$}}\vspace*{-1mm}  \\
\mbox{\scriptsize{$P_1$}}\end{array}\right)}{\longrightarrow}
\begin{array}{c} \cV \\ \oplus \\
\cV \end{array} \stackrel{(P_1~-P_0)}{\longrightarrow} \cV \to 0
\end{equation} 
is forced to be exact (and is split) by the identity $ id_\cV=
Q_1P_0+Q_0P_1$.  The splitting sequence takes the same form with $Q_0$ and
$-Q_1$ formally replacing, respectively, $P_0$ and $P_1$.
For example the system 
$$
  P_0 u=f_0,\ \ P_1 u=f_1 
$$ has the exact integrability condition $P_1f_0 = P_0 f_1$ and if
this holds then the solution $u = Q_1f_0+Q_0f_1$ is unique.  Including the
projections for the bundle in the centre of the sequence \nn{seq1} we
obtain a diagram
$$
\xymatrix@R-25pt{
  &*=0{} \ar@{-}[rr] && *=0{} \ar[dd] \\
 &&\cV   \ar[rd]^-{P_1} \\
   0     \ar[r]
  &\cV   \ar@{-}[uu]
         \ar[ru]^-{P_0}
         \ar[rd]_-{P_1}
 && \cV  \ar[r]
  & 0.   \\
 && \cV  \ar[ru]_-{-P_0}
}
$$
Here the long arrow indicates the
composition $P_1P_0$ and note that by viewing the centre column as a
direct sum we include the information of the original complex.

A related observation follows. This  concerns how, for the very simple
case of $\ell=1$, the essential content of Theorem \ref{inhg} is
captured in the short exact sequence \nn{seq1}. Notation is as above.
\begin{lemma}\label{seq1lem}
If \nn{seq1} is exact and $(u^0,u^1)$ solves the system $P_0 u^0 =f$
and $P_1 u^1=f$, then $(u^1,u^0)=(P_0 u, P_1 u)$ for a unique $u\in
\cV$ satisfying $Pu=f$.
\end{lemma}
\begin{proof}
Note that $(u^1,u^0)$ is in the
null space of $(P_1~-P_0)$ and so, using that \nn{seq1} is exact, we
have  the result.
\end{proof}
\noindent A key point is that this holds without explicit mention of
the splitting $Q$-operators. Of course, for example, $(P_1~-P_0)$ has
a left inverse but we do not any commutativity properties of this
beyond what is forced by \nn{seq1} being exact.  We will return to
this point in Section \ref{compare}.

The case $\ell=2$, i.e.\ the system
\begin{equation} \label{threeEq}
  P_0 u=f_0,\ \ P_1 u=f_1,\ \ P_2 u=f_2 
\end{equation}
demonstrates the general situation more accurately. Here
we have $P=P_0P_1P_2$ where the commutators $[P_i,P_j]$ are all
trivial.  If there is a solution to \nn{threeEq} then it is necessary that $P_if_j = P_jf_i$
for all $0 \leq i,j \leq 2$.
These and further problems with their integrability conditions may be organised into the complex
\begin{equation}\label{3comp}
\xymatrix@R-5pt@C+15pt{
  && \cV \ar[r]^-{P_1} 
         \ar[rd]_-<<*-{\scriptstyle P_2}
   & \cV \ar[rd]^-{P_2} \\
     0 \ar[r]
   & \cV \ar[ru]^-{P_0} \ar[r]_-{P_1} \ar[rd]_-{P_2}
   & \cV \ar[ru]_<<*-{\scriptstyle -P_0} 
         \ar[rd]_<<<*-{\scriptstyle P_2}
   & \cV \ar[r]^-<<<{-P_1}
   &\cV \ar[r]
   & 0 ~.\\
  && \cV \ar[r]_-{-P_1} 
         \ar[ru]_->>*-{\scriptstyle -P_0}
   & \cV \ar[ru]_-{P_0}
}
\end{equation}
Now consider for each operator $\stackrel{\pm P_i}{\longrightarrow}$
also a corresponding operator $\stackrel{\pm Q_i}{\longleftarrow}$ in
the opposite direction. We assume, as before, $[Q_i,P_j]=0$ for
$i,j\in\{0,1,2\}$.  (Note that if also the operators $Q_i$ are
mutually commuting then they too form a complex.)  Denoting the space
of the degree $p \in \{0,1,2,3\}$ by $\cV(p)$ (i.e.\ $\cV(0)=\cV$,
$\cV(1) = \cV \oplus \cV \oplus \cV$ etc.), these complexes become
$$
\xymatrix@C+5pt{
   0      \ar@<2pt>[r]
  &\cV(0) \ar@<2pt>[l]
          \ar@<2pt>[r]^{P(0)}
  &\cV(1) \ar@<2pt>[l]^{Q(1)}
          \ar@<2pt>[r]^{P(1)}
  &\cV(2) \ar@<2pt>[l]^{Q(2)}
          \ar@<2pt>[r]^{P(2)}
  &\cV(3) \ar@<2pt>[l]^{Q(3)}
          \ar@<2pt>[r]
  &0      \ar@<2pt>[l]
}
$$
where the operators $P(p): \cV(p) \to \cV(p+1)$ are given by the corresponding
sum of operators $\pm P_i$ (and similarly for $Q(p)$).
Now we can write the system \nn{threeEq} simply as 
$P(0) u = \mbox{\boldmath$f$}$  where $\mbox{\boldmath$f$} = (f_0,f_1,f_2) \in \cV(1)$.

The cohomology of this complex is related to the solution spaces of
the problems $P(p)u=0$. In particular $H^0 = \cN(P(0))$. 
While in general little could be said about the cohomology the key
point is this.  Observe that 
$$ 
Q(p+1) P(p) + P(p-1) Q(p) = Q_0P_0 + Q_1P_1 + Q_2P_2 
$$ on $\cV(p)$ for $p \in \{0,1,2,3\}$ (where $P(-1)$, $P(3)$, $Q(0)$
and $Q(4)$ are indicated trivial mappings).  Hence if the right hand
side of the last display is equal to the identity then $P(p-1) Q(p)$
is the identity on $\cN( P(p))$ and so the complex is exact; in the 
case that the $P(p)$-complex is exact we shall say the complex is 
a {\em resolution} of the operator
$P$. 
When we have such a resolution then, for example,
the problem $P(0) u = \mbox{\boldmath$f$}$, has a solution only if we
have the ``integrability condition'' $\mbox{\boldmath$f$} \in
\cN(P(1))$ and if this holds then the solution is unique. 

Note that the diagram \nn{3comp} is essentially the Hasse diagram (or
lattice diagram) for the natural poset structure of the power set
$2^{L}$ in the case $L=\{0,1,2\}$.  The situation for a general $\ell
\in \N$ is analogous, and we shall exploit the connection to poset
structure to organise the notation.  So we consider operators $P_i:
\cV \to \cV$, $i \in L := \{0,1,\ldots,\ell\}$ which are mutually
commuting, but otherwise arbitrary. The complex will be constructed
using $2^{|L|}$ copies of $\cV$ as follows. The copies of $\cV$ will
be indexed by subsets $J \subseteq L$, i.e.\ $\cV_J := \cV$, and we
define the operators
$$ P_{J,i} := (-1)^{|J<i|} P_i: \cV_J \to \cV_{J \cup \{i\}}, \quad
   J \subseteq L, i \in L \setminus J $$
where
$$ |J<i| := |\{ j \in J \mid j<i \}|, \quad J \subseteq L. $$
Further we put
$$ \cV(p) := \bigoplus_{\substack{J \subseteq L, \\ |J|=p}} \cV_J
   \quad \mbox{and} \quad
   P(p) := \mbox{\LARGE\bf +}_{p=|J|, i \not\in J} P_{J,i}: \cV(p) \to \cV(p+1). $$

\begin{proposition} \label{exact}
The operators $P(p)$, $0 \leq p \leq \ell$ form a complex
\begin{equation}\label{res}
\xymatrix@C-3pt@M-1pt@R-6pt{
   0           \ar[r]
  &\cV(0)      \ar[r]^-{P(0)}
  &\ \cdots \  \ar[r]^-{P(\ell)}
  &\cV(\ell+1) \ar[r]
  &0.}
\end{equation}
Moreover, if $id_{\cV} = Q_0 P_0 + \ldots + Q_\ell P_\ell$ for some
operators $Q_i: \cV \to \cV$, satisfying $[Q_i,P_j]=0$ for
$i,j=0,\ldots,\ell$, then this complex is exact.
\end{proposition}
\noindent If the complex \nn{res} is exact we shall call it a
resolution of the operator $P$. In the treatment of the Taylor
spectrum for commuting operators on a Banach space the Koszul complex
here is said to be {\em Taylor regular} if it is exact. The main part
of the Proposition here is in Proposition 3, Chapter IV of
\cite{M}. We include the proof here to keep the treatment
self-contained and in terms of a single notational system.

\begin{proof}
First we need to show that $P(p+1) \circ P(p)=0$. This map is a sum of
mappings $P_{J,i,j}: \cV_{J} \to \cV_{J \cup \{i,j\}}$ such that
$|J|=p$ and $i,j \in L \setminus J$ given by the restriction of
$P(p+1) \circ P(p)$ to the source subspace $\cV_J \subseteq \cV(p)$
and projection onto the target subspace $\cV_{J \cup \{i,j\}}
\subseteq \cV(p+2)$.  Fix such a triple $(J,i,j)$ and assume
$i<j$. Denoting $q_i = |J<i|$ and $q_j = |J<j|$, we obtain that
$P_{J,i,j}$ is the sum of the two composite operators $\cV_J\to \cV_{J
\cup \{i,j\}}$ in the following diamond:
\begin{equation}\label{diamond}
\xymatrix@C-3pt@M-1pt@R-6pt{
  & \cV_{J \cup \{i\}}  \ar[rd]^-{(-1)^{q_j+1} P_j} \\
    \cV_J               \ar[ru]^-{(-1)^{q_i} P_i}
                        \ar[rd]_-{(-1)^{q_j} P_j}
 && \cV_{J \cup \{i,j\}}. \\
  & \cV_{J \cup \{j\}} \ar[ru]_-{(-1)^{q_i} P_i} }
\end{equation}
But from this we see immediately that $P_{J,i,j}=0$.

Now assume we have operators $Q_i:\cV\to \cV$ so that $\id_{\cV} = Q_0
P_0 + \ldots + Q_\ell P_\ell$ and $[Q_i,P_j]=0$ as in the
Proposition. Consider for every operator $\stackrel{\pm
P_i}{\longrightarrow}$ also the operator $\stackrel{\pm
Q_i}{\longleftarrow}$ in the opposite direction.  (Then the operators
labelled by $\pm Q_i$ also form a complex, provided $[Q_i,Q_j]=0$ but
we will not need this fact.) We obtain the diagram
$$
\def\dli{\ar@{} |{. \hspace*{2pt} . \hspace*{2pt} . \hspace*{2pt} .}}
\xymatrix@C+0pt{
   0      \ar@<2pt>[r]
  &\cV(0) \ar@<2pt>[l]
          \ar@<2pt>[r]^{P(0)}
  &\cV(1) \ar@<2pt>[l]^{Q(1)}
          \ar@<2pt>[r]^(0.65){P(1)}
  &       \ar@<2pt>[l]^(0.35){Q(2)}
          \hspace{-2ex} \dli[r] \hspace{-2ex}
  &       \ar@<2pt>[r]^(0.4){P(\ell-1)}
  &\cV(\ell) \ar@<2pt>[r]^(0.45){P(\ell)}
          \ar@<2pt>[l]^(0.55){Q(\ell)}
  &\cV(\ell \!+\! 1) 
          \ar@<2pt>[l]^(0.55){Q(\ell+1)}
          \ar@<2pt>[r]
  &0      \ar@<2pt>[l]
}
$$ where $Q(j)$ is the sum of the $\stackrel{\pm Q_i}{\longleftarrow}$
between the corresponding subspaces of $\cV(j)$ and $\cV(j-1)$ for $j
= 1,\ldots,\ell+1$. We denote by $P(-1)$, $P(\ell+1)$, $Q(0)$ and
$Q(\ell+2)$, in an obvious way, the trivial operators at the left and
right extremes of the diagram.  Let us fix $p \in \{0,\ldots,\ell+1\}$
and consider the restriction $Q(p+1) \circ P(p)|_{\cV_J}$ for some $J\subseteq L$ with $|J|=p$.
By definition, $P(p)|\cV_J$ is the sum of
operators $P_{J,i}$ for $i \not\in J$. The $Q$--operators from $\cV_{J
\cup\{i\}} \subseteq \cV(p+1)$ (i.e.\ the target space of $P_{J,i}$)
back to $\cV(p)$ correspond to $j \in J \cup \{i\}$ and have $\cV_{(J
\cup \{i\}) \setminus \{j\}} \subseteq \cV(p)$ as the target space.
For a given $i \not\in J$, the choice $j:=i$ yields the composition
$Q_iP_i: \cV_J \to \cV_J$, and the choices $j \in J$ yield the
operators
\begin{eqnarray} \label{q_ij}
\begin{split}
R_{ij} &= 
  \xymatrix@C+11ex{
     \cV_J  \ar[r]^-{(-1)^{|J<i|}P_i}
    &\cV_{J \cup \{i\}}
            \ar[r]^-{(-1)^{|(J \cup \{i\}) \setminus \{j\} <j|}Q_j}
    &\cV_{(J \cup \{i\}) \setminus \{j\}}
  } \hspace{-8ex} \\
&= q_{ij} P_iQ_j: \cV_J \to \cV_{(J \cup \{i\}) \setminus \{j\}},
   \quad i \not\in J, j \in J.
\end{split}
\end{eqnarray}
where $q_{i,j} \in \{+1,-1\}$ is determined by the previous display. 
Summarising, we have obtained
$$ Q(p+1) \circ P(p)|_{\cV_J} = 
   \Bigl( \sum_{i \not\in J} Q_iP_i \Bigr)_{\mbox{on} \, \cV_J}
   + \sum_{i \not\in J, j \in J}R_{ij}. $$
The same analysis of $P(p-1) \circ Q(p)$ yields
$$ P(p-1) \circ Q(p)|_{\cV_J} =
   \Bigl( \sum_{i \in J} P_iQ_i \Bigr)_{\mbox{on} \, \cV_J}
   + \sum_{i \not\in J, j \in J} R'_{ij} $$
where
\begin{eqnarray} \label{q'_ij}
\begin{split}
R'_{ij} &=
  \xymatrix@C+11ex{
     \cV_J  \ar[r]^-{(-1)^{| J \setminus \{j\} <j|}Q_j}
    &\cV_{J \setminus \{j\}}
            \ar[r]^-{(-1)^{|J \setminus \{j\} <i|}P_i}
    &\cV_{(J \setminus \{j\}) \cup \{i\}}
  } \hspace{-8ex} \\
&= q'_{ij} Q_jP_i: \cV_J \to \cV_{(J \setminus \{j\}) \cup \{i\}},
   \quad j \in J, i \not\in J.
\end{split}
\end{eqnarray}
Summarising again (and using $[P_i,Q_j]=0$), we obtain
\begin{align*}
  &Q(p+1) \circ P(p) + P(p-1) \circ Q(p)|_{\cV_J} = \\
  &=\Bigl( \sum_{i=0}^{\ell} Q_iP_i \Bigr)_{\mbox{on} \, \cV_J} 
   + \sum_{i \not\in J, j \in J} (R_{ij} + R'_{ij}) = \\
  &=\Bigl( \sum_{i=0}^{\ell} Q_iP_i \Bigr)_{\mbox{on} \, \cV_J}
   + \sum_{i \not\in J, j \in J} (q_{ij} + q'_{ij})
     (P_iQ_j)_{\cV_J \to \cV_{(J \setminus \{j\}) \cup \{i\}}}.
\end{align*}
The first sum is the identity according to the assumption. To compute the 
second one we use the explicit form of $q_{ij}$ and  $q'_{ij}$ given by
respectively \nn{q_ij} and \nn{q'_ij}. If $i>j$ then
$$ (-1)^{|J \setminus \{j\} <j|} = (-1)^{|(J \cup\{i\}) \setminus \{j\} <j|}
   \quad \mbox{and} \quad
   (-1)^{|J<i|} = -(-1)^{|J \setminus \{j\} <i|} $$
hence $q_{i,j} = -q'_{i,j}$. One easily sees the latter is true also for
$i<j$. Therefore we obtain 
$$ Q(p+1) \circ P(p) + P(p-1) \circ Q(p) = \id_{\cV(p)}, $$ whence
$P(p-1) \circ Q(p)$ is the identity on $\cN(P(p))$ and the Proposition
follows.
\end{proof}

Note that the identity $\id_{\cV} = \sum_{i=0}^\ell Q_iP_i$ is in
general far weaker than \nn{iddec} required for a decomposition.  This
motivates a rather broader notion of decomposition that we now
introduce.

\subsection{General case: $\al$-decompositions}\label{aldec}
We define and discuss here a generalisation of the notion of a
decomposition which has the decomposition from Section \ref{gendecom}
as simply an extreme (but important) special class.  
Consider the operator $P = P_0 \cdots P_\ell$ from \nn{compP}
and the power set $2^L$ of the index set $L :=
\{0,1,\ldots,\ell\}$. We shall use the notation $P_J := \prod_{j \in
J} P_j$ for $\emptyset \not=J \subsetneq L$ 
and set $P_\emptyset :=
\id_{\cV}$.  Now choose a nonempty subset $\al \subseteq 2^L$ and
assume there exist operators $Q_J: \cV \to \cV$, $J \in \al$ that give
a decomposition of the identity
\begin{equation} \label{al-iddec}
  id_{\cV} = \sum_{J \in \al} Q_J P^J 
\end{equation}
 where
$P^J = P_{L \setminus J}$, and $P_i$s and $Q_J$s satisfy
\begin{equation} \label{al-commut}
  P_iP_j = P_jP_i \quad \mbox{and} \quad 
  P_iQ_J = Q_JP_i \quad i \in L, J \in \al.
\end{equation}

\begin{definition*} For a linear operator $P:\cV\to \cV$, an
expression of the form \nn{compP} will be said to be a
{\em $\al$--decomposition}
of $P$ if the identity \nn{al-iddec} holds with \nn{al-commut} satisfied 
and $\emptyset \not= \al \subseteq 2^L$, $L \not\in \al$.
\end{definition*}

The case of \nn{compP} being a decomposition is a special case of an
$\al$--decomposition with $\al = \{J \subseteq L \mid |J|=1
\}$. 
Toward understanding $\al$-decompositions we employ a 
dual notion of a decomposition, as follows.
\begin{definition*}
We say that $P = P_0 \cdots P_\ell$ is the
\idx{dual $\be$--decomposition}, $\emptyset \not= \be \subseteq 2^L$, 
$\{\emptyset\} \not= \be$ if for 
every $J \in \be$ there exist operators $Q_{J,j} \in \End(\cV)$, $j \in J$ 
such that 
\begin{equation} \label{iddecdual}
  \id_{\cV} = \sum_{j \in J} Q_{J,j} P_j, \quad
  [P_i,P_k] = [Q_{J,j},P_i] =0, \quad
  i,k \in L, j \in J.
\end{equation}
\end{definition*}

Each system $\al \subseteq 2^L$ is partially ordered be restricting
the poset structure of $2^L$. The sets of minimal and maximal elements
in $\al$ will be denoted by $\Min(\al)$ and $\Max(\al)$, respectively.
We say the system $\be \subseteq 2^L$ is a \idx{lower set}, if it is
closed under taking a subset. (That is, if $I \in \be$ and $J
\subseteq I$ then $J \in \be$.) The \idx{upper set} is defined
dually. The lower set and upper set generated by a system $\al
\subseteq 2^L$ will be denoted by 
$\cL(\al) := \{ J\subseteq I~|~ I\in \al\}$ 
and $\cU(\al) :=  \{ J \supseteq I~|~ J \subseteq L~\mbox{and}~I\in \al\}$,
respectively.

\begin{lemma} \label{MinMax}
Let $\al \subseteq 2^L$. Then 
$P = P_0 \cdots P_\ell$ satisfies the following: \\
(i) it is an $\al$--decomposition $\Longleftrightarrow$
it is a $\Max(\al)$--decomposition $\Longleftrightarrow$
it is an $\cL(\al)$--decomposition \\
(ii) it is a dual $\al$--decomposition $\Longleftrightarrow$
it is a dual $\Min(\al)$--decomposition $\Longleftrightarrow$
it is a dual $\cU(\al)$--decomposition.
\end{lemma}

\begin{proof}
The proof of (i) follows easily from the definitions, the proof of (ii)
is also obvious.
\end{proof}

To formulate the relation between $\al$-- and dual
$\al$--decompositions, we need the following notation. We put $\al^u
:= 2^L \setminus \cL(\al)$ 
and $\al^l := 2^L
\setminus \cU(\al)$. Clearly $(\al^u)^l = \cL(\al)$ and $(\al^l)^u =
\cU(\al)$.  
Also it is easily seen that
\begin{eqnarray} \label{UL}
\begin{split}
  &\al^u = \{J \subseteq L \mid \forall I \in \al: 
            J \setminus I \not= \emptyset \} \\
  &\al^l = \{J \subseteq L \mid \forall I \in \al:
            I \setminus J \not= \emptyset \}.
\end{split}
\end{eqnarray}

The first part of the following proposition describes the duality in 
the special case \nn{iddec}. 
\begin{proposition} \label{Peucgen}
\nn{iddec} is equivalent to
\begin{equation} \label{paireucgen} 
  \id_{\cV} = Q_{i,j} P_i +Q_{j,i} P_j  
\end{equation}
where $Q_{i.j}\in \End (\cV)$ and satisfy $[Q_{i,j},P_k]=0$ for every
triple of integers $(i,j,k)$ such that $0 \leq i , j,k \leq \ell$ and
$i \not= j$. That is, \nn{compP} with \nn{iddec} 
is equivalent to the dual $\be$--decomposition for
$\be = \{J \subseteq L \mid |J|=2 \}$.

More generally, \nn{compP} is an $\al$--decomposition if and only if 
it is a dual $\al^u$--de\-composition. Equivalently, \nn{compP} is a
dual $\be$--decomposition if and only if it is a $\be^l$--decomposition.
\end{proposition}

\begin{proof}
We shall prove the first part of the general statement, i.e.\
that \nn{compP} is an $\al$--decomposition if and only if
it is a dual $\al^u$--decomposition. 
Also we will suppose
$\al = \cL(\al)$. This is no loss of generality due to Lemma \ref{MinMax}.

Assume \nn{compP} is 
$\al$--decomposition and consider $J \in \al^u$. That is, 
$J \setminus I \not= \emptyset$ for all $I \in \al$.
For any $I \subseteq L$, we have $L \setminus I = J' \cup (J \setminus I)$,
$J' \subseteq L$ as the disjoint union. Hence
$$ Q_I P^I = Q_I P_{L \setminus I} = (Q_I P_{J'}) P_{J \setminus I}. $$
From this, it is obvious that the identity
$\id_{\cV} = \sum_{I \in \al} Q_I P^I$ can be easily rewritten to the form
\nn{iddecdual} because $\emptyset \not= J \setminus I \subseteq J$
for all $I \in \al$. (The commutation relations in \nn{iddecdual} are clearly
satisfied.)

Now assume \nn{compP} is the dual $\al^u$--decomposition, i.e.\ 
$\id_{\cV} = \sum_{j \in J} Q_{J,j}P_j$ for every $J \in \al^u$.
We shall prove 
that for every $J \in \al^u$
we have a decomposition of the identity
\begin{equation} \label{induction}
  id_\cV = \sum_{I \in \al} Q_{I,J} P_{J \setminus I}, \quad 
  Q_{I,J} \in \End(\cV),  
\end{equation}
such that $[Q_{I,J},P_k]=0$ for every $I \in \al$, and $k \in L$. Then
the proposition follows from the choice $J := L \in \al^u$ as $P_{L
\setminus I} = P^I$. The proof will use induction on the
partial ordering of $2^L$ (given by inclusion), will use that $2^L$ is
the disjoint union $2^L= \al^u \cup \al$, and also that, since
$\al=\cL(\al)$, we have the least element $\emptyset \in \al$.

Before
we do the induction let us first consider as easy case 
which indicates how the argument works, viz.\
 $J \in \Min(\al^u)$. It
 follows from this minimality that for every $j \in J$ we obtain $I_j
 := J \setminus \{j\} \in \al$ hence $\{j\} = J \setminus I_j$ for
 some $I_j \in \al$. Using this and since \nn{compP} is a dual
 $\al^u$--decomposition and $J \in \al^u$, we conclude $\id_\cV =
 \sum_{j \in J} Q_{J,j} P_j = \sum_{j \in J} Q_{J,j} P_{J \setminus
 I_j}$.  But the latter sum is of the form \nn{induction} because $I_j
 \in \al$.  (We put $Q_{I,J} := 0$ for every $I \in \al$ not of the
 form $I_j$ for some $j \in J$.) The commutativity conditions in
 \nn{induction} follow from the definition of the dual
 $\al^u$--decomposition.

Now consider $J \in \al^u$. Since $2^L= \al^u \cup \al$ is a disjoint
union, there are sets $J'$ and $J''$ so that $J = J' \cup J''$ where
$J \setminus \{j\} =: I_j \in \al$ for $j \in J'$, and $J \setminus
\{j\} =: J_j \in \al^u$ for $j \in J''$. Now, as $J \in \al^u$, the assumption 
of a dual $\al^u$-decomposition gives
the identity
$$ id_\cV = \sum_{j \in J'} Q_{J,j} P_{J \setminus I_j} + \sum_{j \in
          J''} Q_{J,j} P_j ~,
$$ 
where we have used that $J \setminus I_j = \{j\}$ for $j \in J'$. This
first sum is of the form required in \nn{induction}, the second one is
not. But since $\al^u \ni J_j \subsetneq J$ for $j \in J''$, we may
assume, by the induction, that $id_\cV = \sum_{I \in \al} Q_{J_j,I}
P_{J_j \setminus I}$.  Acting on this by $P_j$, we obtain $P_j =
\sum_{I \in \al} Q_{J_j,I} P_{J \setminus (I \setminus \{j\}})$.  Here
$I \setminus \{j\} \in \al$ (because $\al = \cL(\al)$) hence the
latter sum is of the form on the right hand side of \nn{induction}.
Consequently, putting these expressions for $P_j$, $j \in J''$ into
the previous display, we obtain decomposition of identity of the form
of \nn{induction}.  The required commutativity relations clearly hold
thus the proposition follows.
\end{proof}
\begin{remark} A main point of the Proposition above is to shed light on 
the nature of $\al$-decompositions. The $\al$-decomposition is what
gets directly used in studying the solution space for $P$. However at
first this seems rather mysterious since, for example, the $P^J$ in
the identity \nn{al-iddec} are complementary to the $P_J$. The first
part of the Proposition exposes one view of what it means to say that
$P_0P_1\cdots P_\ell$ is a decomposition: it shows that \nn{iddec} is
equivalent to the $P_i$s being mutually relatively invertible.
We will see in Section \ref{compare} that this picture generalises.

Next note that the proof
above, begining with \nn{iddecdual}, inductively constructs explicit
formulae for the $Q_J$ in \nn{al-iddec} in terms of products of the
$Q_{J,j}$ from \nn{iddecdual}.  Note also that although, we do not
require $[Q_{J,j},Q_{J',j'}]=0$ for $J,J' \in \al$, $j \in J$, $j' \in
J'$ in \nn{iddecdual}, in the special case \nn{paireucgen} in the
Proposition one shows that $[Q_{i,j},Q_{j,i}]=0$ easily follows from
\nn{paireucgen} and the vanishing of the $[Q_{i,j},P_k]$ as assumed.
\end{remark}

The subsets $\al \subseteq 2^L$ are partially ordered by inclusion
(i.e.\ now we use the poset structure of $2^{2^L}$).  Given an
operator $P$ in the form \nn{compP} consider the family $\Gamma$ of
systems $\al$ such that \nn{compP} is a dual
$\al$--decomposition. Then $\Gamma$ has the greatest element $\al_P =
\bigcup_{\al \in \Gamma} \al$. Then an ``optimal'' choice for the
(dual) $\al$--decomposition of $P$ is $\al := \Min(\al_P)$.  (We want
to have in $\al$ to the smallest possible subsets of $L$. So if the
$P_i$s are not invertible then the case of a dual decomposition may be
regarded as the best we can do. With this philosophy we thus take
$\al_P$.  Then using Lemma \ref{MinMax} we take $\al := \Min(\al_P)$
as it is easier to work with a smaller number of subsets.)
Consequently, we obtain the optimal choice $\be := \Max((\al_P)^l)$
for the $\be$--decomposition of $P$.

In the case one is able to decide, given a subset $J \subseteq L$, whether 
$id_V = \sum_{j \in J} Q_{J,j} P_j$ for some $Q_{J,j} \in \End(\cV)$, it is 
easy to find the optimal (dual) decompositions. This is, for example, the case
of polynomial operators discussed in Section \ref{multivar}.

\vspace{2mm}
 
The decomposition used in Lemma \ref{projlem} and Theorem \ref{inhg}
is a special case of the $\al$--decomposition, $\emptyset \not = \al
\subseteq 2^L$ where $L = \{0,1,\ldots,\ell\}$. The null
spaces of the $P_J$ will in general meet non-trivially. However note
the following.
\begin{lemma}
If $\al\subseteq 2^L$ gives an $\al$-decomposition of $P$
then 
\begin{equation}\label{nproj}
Q_IP^I:\cN(P)\to \cN(P_I) \quad \mbox{ for all } I\in \al~.
\end{equation}
If $\al$ satisfies $I\cap J=\emptyset$ for all $I\neq J\in \al$ then,
for each $I\in \al$, $Q_IP^I $ in \nn{nproj} is a projection.
\end{lemma}
\begin{proof} 
The point is that if the sets in $\al$ are mutually disjoint then $P^J$ (and
hence $Q_JP^J$) annihilates $\cN(P_I)$ whenever $I\neq J$. So the
proof of Lemma \ref{projlem} generalises easily.
\end{proof}

Using the Lemma and by an easy adaption of the proof of Theorem
\ref{inhg} we obtain the following.

\begin{theorem}\label{winhg}
Assume $P:\cV\to \cV$ as in \nn{compP} is an $\al$--decomposition.  
Let us fix $f\in \cV$.  
There is a surjective mapping $B$ from the space of solutions 
$(u_J)_{J \in \al} \in \oplus^{|\al|} \cV$ of the
problem
\begin{equation}\label{wmultiprobg} 
 P_J u_J =f, \quad J \in \al. 
\end{equation}
onto the space of solutions $u\in \cV$ of $P u=f$.

Writing $\cV_P^f$ for the solution space of $P u=f$ and (for $J \in \al$) 
$\cV^f_J$ for the
solution space of $P_J \tilde{u}=f$.  The map
$B:\times_{J \in \al}  \cV^f_J \to \cV^f_P$ is given by
$$ (u_J)_{J \in \al} \mapsto \sum_{J \in \al} 
   Q_J  u_J ~.
$$
A right inverse for this is 
$F:\cV^f_P \to \times_{J \in \al} \cV^f_J$  given (component-wise) by
$$ u\mapsto P^J u  ~;$$
on $\cV$ we have $B\circ F=id_{\cV}$. 

If $\al$ satisfies $I\cap J=\emptyset$ for all $I\neq J\in \al$ then,
$F$ is a 1-1 mapping and $F\circ B$ is the identity on the solution space to  
\nn{wmultiprobg}.
\end{theorem}

\def\bigplus{\mbox{\large$+$}}

\noindent Hence the  generalisation of Corollary \ref{maincor} is as follows.
\begin{corollary}\label{wmaincor}
For $P:\cV\to \cV$, with \nn{compP} giving an $\al$--decomposition, 
$\emptyset \not = \al \subseteq 2^L$, we have
$$ \cR(P)=\bigcap_{J \in \al} \cR(P_J), \quad
   \cN(P)=\bigplus_{J \in \al} \cN(P_J). $$
If $\al$ consists of mutually disjoint sets then we have
$$
\cN(P)=\bigoplus_{J \in \al} \cN(P_J),
$$
and this is given by 
$$
id_{\cN(P)}=\sum_{J\in \al} \Proj_J
$$ where, for each $I\in\al$, $\Proj_I:\cN(P)\to
\cN(P_I)$ is the projection given by the restriction of $Q_IP^I$ from
\nn{al-iddec}.
\end{corollary}

\noindent So although the assumption of an $\al$--decomposition, for
an operator $P$, is in general a vastly weaker requirement than that of
a decomposition, we still have the critical result that one may solve the
inhomogeneous problem $P u=f $ by treating a ``lower order'' problem
involving the same inhomogeneous term $f$.

\vspace{1ex}


\subsection{$\al$-decompositions in terms of operator resolutions}
\label{compare}
Recall that the complex \nn{res} in Proposition \ref{exact} promotes
to being an operator resolution (i.e.\ is exact) if we make the the
rather weak assumption $id_\cV=\sum_0^\ell Q_i P_i$ (with the usual
commutativity of operators assumed). On the other hand expression
\nn{paireucgen} in Proposition \ref{Peucgen} shows when one has a
decomposition $P=P_0P_1\cdots P_\ell$ (i.e.\ \nn{iddec} holds) then
every diamond subcomplex \nn{diamond} of the operator resolution
diagram \nn{res} is exact. In a sense, that we now make precise, this
is the key algebraic content of a decomposition.

Consider then $P=P_0P_1\cdots P_\ell$, where as usual the $P_i\in \End
\cV$ are mutually commuting. We have the complex \nn{res}. Let us
assume that in this each diamond subcomplex of the form \nn{diamond}
is exact (in the sense of \nn{seq1}). Then the complex \nn{res} is
exact and so gives an operator resolution. We shall investigate to
what extent the results for decompositions survive if we take this
setting without explicitly requiring the identity \nn{iddec}.

 We earlier discussed the case $\ell=1$. To shed light on the general
situation we look now at the case $\ell=2$, $P=P_0P_1P_2$, so we have
the complex \nn{3comp}.  For $f\in \cV$, consider the inhomogeneous
problem $P_0v_0=f$, $P_1v_1=f$, $P_2v_2=f$. Since each diamond is
exact we have that
$$
\left(\begin{array}{c} v_0\\
v_1
\end{array}\right)= \left(\begin{array}{c} -P_1 u_2\\
-P_0 u_2
\end{array}\right) \quad 
\left(\begin{array}{c} v_1\\
v_2
\end{array}\right)= \left(\begin{array}{c} P_2 u_0\\
P_1 u_0
\end{array}\right)~.
$$ These are consistent only if $P_2 u_0+P_0u_2=0$, and, when this holds, using
that the $P_0$, $P_2$ diamond is exact we find that, for $i=0,1,2$,
$v_i=P^i u$ for $u\in \cV$ satisfying $P u=f$ (cf.\ Theorem
\ref{inhg}).  

The results for $\ell=1,2$ extend to general $\ell\in
\mathbb{N}$.
 
\begin{theorem}\label{Qfree}
Suppose that we have $P=P_0P_1\cdots P_\ell$, as in \nn{compP}.
Suppose also
that in the corresponding sequence \nn{res} every diamond \nn{diamond}
is exact in the sense of \nn{seq1}. Then all results of Theorem
\ref{inhg} hold except the map $B$ should replaced by the map $B'$ given
$(u^0,\cdots ,u^\ell)\mapsto u$  by taking in \nn{res} the unique
preimage (of the map $F$) in $\cV_\emptyset$ of
$$ (u^0,\cdots , u^\ell)\in \cV_{L\setminus \{0\}}\oplus \cdots \oplus
\cV_{L\setminus \{\ell\}}
$$
solving 
\begin{equation}\label{star}
P_i u^i=f, \quad\quad i=0,1,\cdots ,\ell.
\end{equation}

We have
$$
\cR(P)=\bigcap^{i=\ell}_{i=0}\cR(P_i), \quad
\cN(P)\cong\bigoplus^{i=\ell}_{i=0}\cN(P_i).
$$
\end{theorem}

\begin{proof} First note that since each diamond in the sequence
  \nn{res} is a complex then $P_iP_j=P_jP_i$ for $i,j\in \{0,\cdots
  \ell \}$ and the sequence is a complex.

  If $f\in \cR(P)$ and $P u=f$ then recall that $(u^0,\cdots
  ,u^\ell):=(P^0 u,\cdots, P^\ell u)$ is a solution of \nn{star}. We
  will (strong induction to) prove that any solution of \nn{star} has
  this form, as forced by the consistency of exact diagram \nn{res}.
  Note that by Lemma \ref{seq1lem} this is true for case $\ell=1$.

Assume now that $\ell\geq 2$.  Starting at $\cV_{\{ 0\}}$ and
$\cV_{\{1 \}}$, in the (length $\ell+1$) resolution diagram \nn{res}
for $P_0\cdots P_\ell$, there are subcomplexes of length $\ell$ that
each take the form of \nn{res}; in both of these the terminal space is
$\cV_L$ (where, as usual, $L:=\{0,1,\cdots ,\ell\}$). By the inductive 
hypothesis, consistency of these
subcomplexes mean that there is $u_0\in \cV_{\{ 0\}}$ satisfying $P^0
u_0=f$ and similarly $P^1 u_1 =f$. 
(Recall $P^i$ means $P/P_i$.) Now $\cV_{\{ 0,1\}}$ is in both
subcomplexes and we obtain a consistency condition: by the process of
repeatedly using Lemma \ref{seq1lem} to take preimages and enforce
consistency (at each diamond) in order to solve for $u_0$ and $u_1$,
it follows easily that $u_0$ and $u_1$ are both ``potentials'' for the
(by induction unique) entry in \ $\cV_{\{0,1\}}$. Hence
$$
\left(\begin{array}{cc}-P_1 &P_0 \end{array}\right)
\left(\begin{array}{c} u_0\\ u_1 \end{array}\right) =0\in \cV_{\{0,1\}}~.
$$ Since the diamond \nn{diamond} for $J=\emptyset$ is exact (in the
sense of \nn{seq1}) it follows that necessarily $u_0=P_0u$ and
$u_1=P_1 u$ for some $u\in \cV_{\emptyset}$.  From these it follows,
respectively, that $u^i=P^i u$, for $i=1,\cdots ,\ell$ and $u^j = P^j
u$, for $j=0,2,3,\cdots ,\ell$. It also follows that $P u=f$.

By construction $B'$ is 1-1 and $F\circ B'$ is the identity on the
solution space.  That the forward map (in the notation of Theorem
\ref{inhg}) $F$ is 1-1 is an easy consequence of the injectivity of
$(P_i,P_j):\cV \to \oplus^2 \cV$ for each pair distinct pair $(i,j)\in
L\times L$. 
The direct sum in last
display follows for the same reason.
\end{proof}

It seems likely that there are analogous simplifications for the
general $\al$-de\-compositions. It has also not escaped our attention
that these ideas suggest that there should be extensions of the ideas
here to the setting where one has a suitable commuting diagram but
without assuming that diagram is constructed from commuting operators.
This will be
taken up elsewhere.  Note that although the Theorem here is
conceptually powerful and a far stronger result overall than Theorem
\ref{inhg}, it seems likely that in practice the identity \nn{iddec}
with \nn{commut} is rather useful. In particular one then obtains the
projections in Lemma \ref{projlem}. Also, as we shall see in the
following sections, for a large class of operators we have have all
these identities algebraically.

The resolution diagrams give us a ``pictorial'' understanding of the
$\alpha$-de\-composit\-ions.  For each dual $\beta$-decomposition
$P=P_0P_1\cdots P_\ell$ and $J\in \beta$ we have $id_\cV=\sum_{j\in
J}Q_{J,j}P_j$ (with appropriate commutativity conditions) and so a
collection of length $|J|$ exact subcomplexes of the resolution for
$P$. Each of these is itself an operator resolution for $\prod_{j\in
J}P_j $. The size of the $|J|$ as we range over $J\in \beta$ gives
some measure of the strength of the dual $\beta$-decomposition: the
smaller the sets $J\in \al$ the stronger the decomposition.  For
example $\beta=\{\{0\}, \{ 1\}, \cdots \{ \ell\} \}$ is the case that
all the $P_i$ are invertible.  The duality in Proposition
\ref{Peucgen} allows us therefore to understand $\al$-decompositions
in the same way: Small sets $I$ in $ \al$ indicate a strong
decomposition.

\begin{remark}
  Note that the complexes \nn{res} discussed in Section
  \ref{complexes} were constructed from an arbitrary set
  $P_0,\ldots,P_\ell$ of mutually commuting endomorphisms of $\cV$.
  Hence using the notation used in \nn{al-iddec}, we can take this set
  to be $\{P^J | J \in \al \}$ for some nonempty system $\al \in 2^L$.
  Then it follows immediately from the above proposition that if
  \nn{compP} is an $\al$--decomposition then the corresponding complex
  is exact.
\end{remark}

\begin{remark}
As a final point we note that there are other approaches to the
inhomogeneous case that na\"\i vely seem similar to Theorem
\ref{inhg}. For example note the following.  Assume $P$ to be in the
form \nn{compP} (with the factors not necessarily commuting). Then
clearly $Pu = f$ has a solution if and only if there is a sequence
$f_0,\ldots,f_{\ell} \in \cV$ satisfying
\begin{equation} \label{sequence}
  P_0 f_0 = f, P_1 f_1 = f_0,
  \ldots, P_\ell f_{\ell} = f_{\ell-1}. 
\end{equation}
 So it is sufficient to find such a sequence to obtain a solution
$u=f_\ell$ of $Pu = f$.  
 However this is simply a variant of the
idea from differential equation theory where, through the introduction
of new variables, one replaces a differential equation by a system of
lower order equations. This is very different from Theorem \ref{inhg}.
The system here does not replace $Pu=f$ with a new inhomogeneous
equation, but rather replaces it with a sequence of problems.  We do
not have the ``source term'' $f_{0}$ in $P_1f_{1}=f_0$ until
we have solved the previous problem $P_0f_{0}=f$ and so on.
 \end{remark}

\section{Algebraic decompositions}\label{algdec}

Here we consider operators $P$ polynomial in mutually commuting
operators $D_0, \cdots, D_k$. In this setting we show that generically
we obtain $\al$-decompositions. In fact in this Section we derive those
decompositions (and $\al$-decompositions) that may be obtained in a
purely algebraic or algebraic-geometric manner from the polynomial
formula for the operator. Thus these are universal results that are
{\em independent of the operators} $D_0, \cdots, D_k$. An important feature
of these cases is that the ``relative inverses'', viz.\ the
$Q$-operators in \nn{iddec} and Theorem \ref{winhg}, are then also
obtained as operators {\em polynomial in the same operators} $D_0,
\cdots, D_k$. Thus if, for example, we dealing with $P$ a differential
operator then these relative inverses, are also differential
operators.

The simplest setting and the strongest
results are obtained in the case of operators polynomial in a single
other operator. Here we derive explicit formulae for the decomposition
that are significantly simpler and more efficient than expected from
the general setup.

\subsection{Operators polynomial in a single operator $\cD$} \label{mainsec}

\noindent Let $\cV$ be a vector space over the field $\bF$.  Suppose
that $\cD:\cV\to \cV$ is a non-trivial linear endomorphism. We may
consider the commutative algebra $\bF[\cD]$ of consisting of those
endomorphisms $\cV\to \cV$ which may be given by expressions
polynomial (with coefficients in $\bF$) in $\cD$. Clearly there is an
algebra epimorphism from $\bF[x]$ onto $\bF[\cD]$ given by mapping a
polynomial $P[x]=\sum_{i=0}^k \kappa_i x^i$ to the operator $P[\cD]$,
a formula for which is given by formally replacing the indeterminate
$x$ in $P[x]$ by $\cD$. That is, a formula for $P[\cD]:\cV\to \cV$ is
$\sum_{i=0}^k \kappa_i \cD^i $ where we write $\cD^i$ as a shorthand
for the $i$-fold composition of $\cD$. This algebra map sends
$1\in \bF[x]$ to $id_{\cV}$.

We begin by treating operators of the form \nn{Popform}. That is
$P=P_0P_1\cdots P_\ell$ where $P_i=(\cD+\la_i)^{p_i}$, with the
$\lambda_i\in \bF$ mutually distinct and for $i=0,\cdots ,\ell$,
$p_i\in {\Bbb Z}_{\geq 1}$. Since the algebra $\bF[\cD]$ is
commutative, we may access the results of Section \ref{gendecom}
provided we obtain the identity \nn{iddec}. This we have
from the Euclidean algorithm as follows.
To a polynomial of the form
\begin{equation}\label{Pform}
P[x]=(x+\lambda_0)^{p_0}(x+\lambda_1)^{p_1} \cdots (x+\lambda_\ell)^{p_\ell}
\end{equation} 
(where the $\lambda_i\in \bF$ are mutually distinct and, for
$i=0,\cdots ,\ell$, $p_i\in {\Bbb Z}_{\geq 0}$) we have the following
decomposition of the unit in $\bF[x]$.  We write
$P_i[x]:=(x+\la_i)^{p_i}$ and then $P^i[x]$ for the polynomial
$P[x]/P_i$, $i=0,1,\cdots ,\ell$.
\begin{lemma}\label{step}
There exist polynomials $Q_i[x]$, each of degree at
most $(p_i-1)$, so that 
$$
1=Q_0[x]P^0[x]+Q_1[x]P^1[x]+\cdots +Q_\ell[x]P^\ell[x].
$$
\end{lemma}
\noindent Note that it is also easy to give a short inductive proof of this. The key
specialisation here is the bound on the degree of the $Q_i$s,
otherwise the display is immediate from the polynomial variant of
Proposition \ref{Peucgen}.

From this Lemma we immediately have specialisations of Theorem \ref{inhg} and
Corollary \ref{maincor}. However before we write these 
we would like explicitly to give formulae for the $Q_i[x]$ in the
Lemma. We derive these in way which is rather suitable to our proposed
applications. 
First observe that if
$\cB$ is an operator on $\cV$ then, for $p\geq 2$, the solution space in $\cV$
of $\cB^p u=0$ 
includes, for example, $u$ such that $\cB u=0$. 
The solution space of $\cB^p u=0$ is filtered.
Given a solution
$u$  we may obviously write $u$ as a sum 
\begin{equation}\label{fexp}
u=u^{(0)}+u^{(1)}+\cdots +u^{(p-1)}
\end{equation} 
where $\cB^{p-s}u^{(s)}=0$, but such expansions are not unique.
For example for any $\alpha \in\bF $ we may take $u^{(0)}=(u-\alpha\cB u)$ and
$u^{(1)}= \alpha\cB u$.

Next  observe that
we may think of $(\cD+\lambda)$ as a nilpotent operator on the
solution space $\cV_{\lambda}$ of \nn{power}. Thus for $\la\neq \mu\in \bF$ 
the operator
$(\cD+\mu)$ is polynomially invertible on this space.  With
\begin{equation}
(\cD+\mu)_\lambda^{-1}:={(\mu-\lambda)^{-1}}\big(1
+\frac{(\cD+\lambda)}{(\la-\mu)} +\cdots
+\frac{(\cD+\lambda)^{(p-1)}}{(\la-\mu)^{(p-1)}}
\big)
\label{invform}
\end{equation}
we have 
$$
(\cD+\mu)_\lambda^{-1} (\cD+\mu)=id_{\cV_\lambda}. 
$$
We shall write $(\cD+\mu)_\lambda^{-2}$ to mean
$(\cD+\mu)_\lambda^{-1}\circ (\cD+\mu)_\lambda^{-1}$ and so forth.

Now we construct the polynomial analogues of
\nn{invform}.  Write $Q[x]$ to denote the polynomial
$(x+\lambda)^p$. Suppose we consider $\bF[x]/\langle Q[x] \rangle$ meaning the
algebra of polynomials modulo the ideal generated by $Q[x]$. As a
multiplication operator on $\bF[x]/\langle Q[x] \rangle$, 
$(x+\lambda)$ is nilpotent and if, with $\mu$ as above, we write
\begin{equation}
(x+\mu)_\lambda^{-1}:={(\mu-\lambda)^{-1}}\big(1
+\frac{(x+\lambda)}{(\la-\mu)} +\cdots
+\frac{(x+\lambda)^{(p-1)}}{(\la-\mu)^{(p-1)}}
\big)
\label{polyinvform}
\end{equation}
then,
$$
(x+\mu)_\lambda^{-1} (x+\mu)=1 \mod \langle Q[x] \rangle. 
$$

Similarly considering $P[x]$ as in \nn{Pform} note that in
$\bF[x]/\langle P[x] \rangle$ the polynomial $(x+\la_i)$ is nilpotent as a
multiplication operator on $\big[\prod_{i\neq j =0}^{j=\ell}
(x+\lambda_j)^{p_j}\big]$, since $(x+\la_i)^{p_i}\big[\prod_{i\neq j
=0}^{j=\ell} (x+\lambda_j)^{p_j}\big]=P[x]$.
Consider now the vector in $\bF[x]/\langle P[x] \rangle$ given by, 
\begin{equation}\label{Prform} 
\Pr_i[x] := \big[ \prod_{i\neq j =0}^{j=\ell}
(x+\lambda_j)^{-p_j}_{\lambda_i}\big] \big[\prod_{i\neq k =0}^{k=\ell}
(x+\lambda_k)^{p_k}\big]~.
\end{equation} 
We will also view this as a multiplication operator on $\bF[x]/\langle
P[x] \rangle$.  Now from Lemma \ref{step} we have $1 =\sum_{i=0}^\ell
Q_i[x]P^i[x]$ in $\bF[x]/\langle P[x] \rangle$, where it should be
noted we use the same notation for the polynomials $1$ and $Q_i[x]$
and so forth as well as for their image in $\bF[x]/\langle P[x]
\rangle$.  Applying $\Pr_i[x]$ to both sides of this identity we have
$$
\Pr_i[x] 1= \Pr_i[x] Q_i[x]P^i[x] \mod \langle P[x] \rangle
$$ since if $k\in \{0,\cdots ,\ell\}$ is distinct from $i$ then
$\Pr_i[x] P^k[x]$ vanishes modulo $\langle P[x] \rangle$.  But now
note that $(x+\la_i)$ is nilpotent on $Q_i[x]P^i[x]$, in $\bF[x]/\langle
P[x] \rangle$, as $(x+\la_i)^{p_i} P^i[x]=P[x]$.
  Thus if $i\neq j\in \{0,\cdots , \ell\}$ then, for
example, $(x+\lambda_j)^{-1}_{\lambda_i} (x+\la_j)$ acts as the
identity on $Q_i[x]P^i[x]$.  Hence
$$
\Pr_i[x] Q_i[x]P^i[x] =Q_i[x]P^i[x] \mod \langle P[x] \rangle,
$$
and so 
$$
\Pr_i[x]=\Pr_i[x] 1=Q_i[x]P^i[x] \mod \langle P[x] \rangle.
$$
Thus from the Lemma and these observations we have 
\begin{equation}\label{modver}
1=\sum_{i=0}^\ell \Pr_i[x] \mod \langle P[x] \rangle.
\end{equation}

Finally we note that we may normalise the formula for
$\Pr_i[x]$. Each term $(x+\lambda_j)^{-1}_{\lambda_i}$ in the product 
$$
\prod_{i\neq j =0}^{j=\ell}
(x+\lambda_j)^{-p_j}_{\lambda_i}
$$ is a sum of powers of $(x+\la_i)$. In
$(x+\lambda_j)^{-1}_{\lambda_i}$ it is only necessary to keep these
powers up to $(x+\la_i)^{p_i-1}$ as, recall,
$(x+\la_i)^{p_i}\big[\prod_{i\neq j =0}^{j=\ell}
(x+\lambda_j)^{p_j}\big]=P[x]$. Similarly, since we are applying the
result to $\big[\prod_{i\neq j =0}^{j=\ell} (x+\lambda_j)^{p_j}\big]$
and calculating modulo $\langle P[x] \rangle$, we may then expand the
product $\prod_{i\neq j =0}^{j=\ell} (x+\lambda_j)^{-p_j}_{\lambda_i}$
writing the result as a linear combination of powers of $(x+\la_i)$
but always keeping only powers $(x+\la_i)^q$ for $q\in {\Bbb Z}_{\geq
0}$ such that $q\leq p_i-1$.  Let us write
$$
N\big(\big[ \prod_{i\neq j =0}^{j=\ell}
(x+\lambda_j)^{-p_j}_{\lambda_i}\big]\big)
$$
for this normalised formula for $Q_i[x]$. Thus we have
$$
\Pr^N_i[x]=N\big(\big[ \prod_{i\neq j =0}^{j=\ell}
(x+\lambda_j)^{-p_j}_{\lambda_i}\big]\big) \big[\prod_{i\neq j =0}^{j=\ell}
(x+\lambda_j)^{p_j}\big]
$$ for the corresponding normalised formula for $\Pr_i[x]$,
$i=0,\cdots,\ell$.  So we have
$$
1= \sum_{i=0}^\ell N\big(\big[ \prod_{i\neq j =0}^{j=\ell}
(x+\lambda_j)^{-p_j}_{\lambda_i}\big] \big)\big[\prod_{i\neq k =0}^{k=\ell}
(x+\lambda_k)^{p_k}\big] \mod \langle P[x] \rangle.
$$ But now observe that the normalised formula in the display has
degree at most $p-1$ in $x$, where $p$ denotes  the degree of $P[x]$.
Thus we have the following result.
\begin{theorem}\label{polyid}In $\bF[x]$ we have the identity 
\begin{equation}\label{genpolyid} 
1=\sum_{i=0}^\ell \Pr^N_i[x]=\sum_{i=0}^\ell N\big(\big[ \prod_{i\neq j =0}^{j=\ell}
(x+\lambda_j)^{-p_j}_{\lambda_i}\big]\big) \big[\prod_{i\neq k =0}^{k=\ell}
(x+\lambda_k)^{p_k}\big]~.
\end{equation}
We associate this to the polynomial \nn{Pform}.
\end{theorem}
\noindent (It seems likely that this identity is known from the theory
of partial fractions.)  It follows that for $P$ of the form
\nn{Popform} we have the decomposition identity \nn{iddec} with $Q_i:=
N\big( \prod_{i\neq j =0}^{j=\ell}
(\cD+\lambda_j)^{-p_j}_{\lambda_i}\big)$ and $P^i:= \prod_{i\neq j
=0}^{j=\ell} (\cD+\la_j)^{p_j}$, as follows.
\begin{corollary}\label{genlinid}
Let $\cV$ be a vector space over a field $\bF$. 
Suppose that $\cD:\cV\to \cV$ is a linear endomorphism and that 
$\la_0,\la_1, \cdots  , \la_\ell\in \bF$ are mutually distinct.  
We have the identity in $\End (\cV)$:
$$
id_{\cV}=\sum_{i=0}^\ell \Pr^N_i[\cD] 
=\sum_{i=0}^\ell N\big(\big[ \prod_{i\neq j =0}^{j=\ell}
(\cD+\lambda_j)^{-p_j}_{\lambda_i}\big]\big) \big[\prod_{i\neq k =0}^{k=\ell}
(\cD+\lambda_k)^{p_k}\big]~
$$ 
where 
for $i=0,\cdots ,\ell$, $p_i\in {\Bbb Z}_{\geq 1}$.
\end{corollary}

We obtain immediately 
the following specialisations of the results from Section \ref{gendecom}. 
\begin{theorem} \label{pf0:gen}
Let $\cV$ be a vector space over a field $\bF$. 
Suppose that $\cD:\cV\to \cV$ is a linear endomorphism and consider 
$P:\cV\to \cV$ given by \nn{Popform} with $\la_0,\ldots,\la_{\ell}$
mutually distinct.
Then there is a canonical and unique direct sum decomposition of the 
the null space for $P$,
\begin{equation}\label{sum:gen}
\cV_P=\oplus_{i=0}^\ell \cV_{\lambda_i}~,
\end{equation}
where, for each $i$ in the sum, $\cV_{\lambda_i}$ is the solution
space for $(\cD+\lambda_i)^{p_i}$. This is executed by a canonical
decomposition of the identity on $\cV_P$ 
$$
id_{\cV_P}=\sum_{i=0}^\ell \Proj_i
$$
where  
$\Proj_i:\cV_P\to \cV_{\lambda_i}$, $i=0,\cdots ,\ell$,   are projections
given by the formula
\begin{equation}\label{projform:gen}
\Proj_i:=\big[ \prod_{i\neq j =0}^{j=\ell} (\cD+\lambda_j)^{-p_j}_{\lambda_i}\big] \big[\prod_{i\neq k =0}^{k=\ell} (\cD+\lambda_k)^{p_k}\big]~.
\end{equation}
\end{theorem}
\noindent Note in the theorem we have used the fact that we may omit
the normalisation of \nn{Prform}, since $P[\cD]$ annihilates $u$.
For  the inhomogeneous problems $P u =f$: 
\begin{theorem}\label{inh} Let $P$ be as above. 
Let us fix $f\in \cV$.  There is a 1-1 relationship 
between 
solutions $u\in \cV$ of $P u=f$ and solutions 
$(u_0,\cdots ,u_\ell)\in\oplus^{\ell+1}\cV$ of the problem 
\begin{equation}\label{multiprob} 
(\cD+\lambda_0)^{p_0}u_0=f, \cdots , (\cD+\lambda_\ell)^{p_\ell}u_\ell =f. 
\end{equation}

Writing $\cV_P^f$ for the solution space of $P u=f$ and (for
$i=0,\cdots ,\ell$) $\cV^f_{\la_i}$ for the solution space of 
$(\cD+\la_i)^{p_i}\tilde{u}=f$.
The map $F:\cV^f_P\to \times_{i=0}^\ell\cV^f_{\la_i} $ is given by 
$$
u\mapsto (\prod_{0\neq j=0}^{j=\ell} (\cD+\la_j)^{p_j}, \cdots ,\prod_{\ell\neq j=0}^{j=\ell} (\cD+\la_j)^{p_j} ) ~,
$$
with inverse $B:\times_{i=0}^\ell\cV^f_{\la_i}\to \cV_P^f$ given by
$$
(u_0,\cdots, u_\ell)\mapsto \sum_{i=0}^{i=\ell} 
N \bigl( \prod_{i\neq j=0}^{j=\ell} (\cD+\la_j)^{-p_j}_{\lambda_i} \bigr) u_i .
$$
\end{theorem}

\begin{remark} 
 Tuning our earlier discussion to the current setting
we could opt to expand each $u_i\in \cV_{\lambda_i}$ with respect to
the canonical filtration, say $u_i=u_i^{(0)}+\cdots
+u_i^{(p_i-1)}$. Although such expansions are not unique, we note here that
the explicit form of the projection $\Proj_j$ given in \nn{projform:gen}
gives such an expansion determined canonically by $P$.  The point is
this. Let us fix $j\neq i$ and write $\cV_{\lambda_i}^{(s)}$
for the subspace of elements vectors $h$ in $\cV_{\lambda_i}$
satisfying $(\cD+\lambda_i)^{p_i-s}h=0$. First $\prod_{i\neq j
=0}^{j=\ell} (\cD+\lambda_j)^{p_j} u$ is in
$\cV_{\lambda_i}=\cV_{\lambda_i}^{(0)}$. Thus from
\nn{invform} it follows that
$$
\begin{aligned}
(\cD&+\lambda_j)^{-1}_{\lambda_i} 
\big[\prod_{i\neq j =0}^{j=\ell} (\cD+\lambda_j)^{p_j}\big]u=\\
&
{(\lambda_j-\lambda_i)^{-1}}\big(1
+\frac{(\cD+\lambda_i)}{(\lambda_i-\lambda_j)} -\cdots
+\frac{(\cD+\lambda_i)^{(p-1)}}{(\lambda_i-\lambda_j)^{(p-1)}}
\big)\big[\prod_{i\neq j =0}^{j=\ell} (\cD+\lambda_j)^{p_j}\big]u
\end{aligned}
$$ has the form
$$
h^{(0)}+\cdots +h^{(p_i-1)}
$$ where $h^{(s)}\in \cV^{(s)}_{\lambda_i}$. Now $(\cD+\lambda_i):
\cV_{\lambda_i}^{(s)}\to\cV_{\lambda_i}^{(s+1)}$ where we view
$s\in\mathbb{Z}_{p_i}$. Thus subsequent applications of
$(\cD+\lambda_k)^{-1}_{\lambda_k} $ ($k\neq i$ and $k\neq j$) 
preserve this form and
yield, in the end, an expression 
$$ u_i = \alpha_i \Bigl(
   1 + \alpha_{i,1} (\cD+\la_i) + \ldots +  \alpha_{i,p_i-1} (\cD+\la_i)^{p_i-1}
   \Bigr)
   \left[ \prod_{i \not= j =0}^{j=\ell} (\cD+\la_j)^{p_j} \right] u
$$
where $\alpha_i$ and $\alpha_{i,j}$ are determined 
explicitly by this process, and in fact it is easily seen that 
$$ \alpha_i = \prod_{i\neq j =0}^{j=\ell}
\frac{1}{(\la_j-\la_i)^{p_j}}. \hspace*{1cm} \endrk$$
 \end{remark}

An important (generic) case of Theorems \ref{pf0:gen} and \ref{inh} is
when we have
$$
P= (\cD+\la_0)(\cD+\la_1)\cdots (\cD+\la_\ell)
$$ 
with the $\la_i$ mutually distinct. Then the situation simplifies as follows.
\begin{proposition}\label{pf0}
Let $P$ be as in \nn{Popform} with $p_0=p_1=\cdots =p_\ell=1$.
Then for $i=0,\cdots ,\ell$ 
$$
Q_i= \prod_{i\neq j =0}^{j=\ell} \frac{1}{\la_j-\la_i}~.
$$
\end{proposition}
\noindent This follows immediately from the discussion in the remark
above or is easily verified directly.

\vspace{3mm}

For a linear operator $P:\cV\to \cV$ let us say that $\mu\in \bF$ is
in the spectrum of $P$ ($\mu\in \Spec P$) if $(P-\mu):\cV\to \cV$ is
not invertible (since we are not assuming that $\cV$ is a Banach
space).  Suppose that $P-\mu$ has the form
$\prod_{i=0}^\ell(D+\lambda_i)$ (with the $\lambda_i\in \bF$ {\em not}
necessarily distinct) for some linear operator $\cD:\cV\to \cV$.
Then, since all factors commute, $(P-\mu)$ is injective (surjective)
if and only if each of the factors $(D+\lambda_i) $ is injective
(resp.\ surjective).  Thus if $\bF$ is an algebraically closed field
and $P=P[\cD]$ is polynomial in $\cD$ then the spectrum of $P$ is
obviously generated by the spectrum of $\cD$; $\mu\in \Spec P$ if and
only if $\mu=P[\lambda]$ where $\lambda\in \Spec \cD$.

From the
Theorem \ref{fundthm}, the  eigenspaces are determined by
the generalised eigenvectors of $\cD$.  We
assume $\cD$ to be a linear endomorphism operator on $\cV$, a vector
space over an algebraically closed field $\bF$, in the summary here.

\begin{corollary} \label{spectralthm}
Let $P=P[\cD]$ be polynomial in $\cD$. Then $(\mu,u)$ is an
eigenvalue, eigenvector pair for $P$ if and only if for some $k\in\{1,\cdots ,{\rm deg}(P)\}$
$$
u=u_1+\cdots+u_k~, \quad \quad 0\neq u_i,\quad i=1,\cdots , k,
$$
where, for each $i\in\{1,\cdots , k\}$, 
$(\cD-\la_i)^{p_i} u_i =0$ and $\la_i$ is a multiplicity $p_i$ solution of 
of the polynomial equation $(P-\mu)[x]=0$.
\end{corollary}
\noindent Of course one could study generalised eigenspaces for $P$ in
the same way.

\vspace{2mm}

\vspace{2mm}

\subsection{The real case} If we work over a field that is not
algebraically closed then the situation, in general, is different from
Theorem \ref{fundthm}, since the polynomial $P[x]$ may not factorise
fully.  However all is not lost. We illustrate the situation in the
case that $\bF$ is $\mathbb{R}$, the field of real numbers. This case
can be dealt with via complexification. By viewing $P[x]$ as a
polynomial in $\mathbb{C}[x]$ from the fundamental theorem of algebra
we obtain a factorisation.
\begin{equation}\label{Popformreal}
P[x]= 
\Big( \Pi_{i=0}^{i=\ell_1} (x+\la_i)^{p_{\la_i}}\Big)
\Big( \Pi_{m=0}^{m=\ell_2} (x+\ka_m)^{p_{\ka_m}}(x+\bar{\ka}_m)^{p_{\ka_m}}\Big)
\end{equation} 
Here the $-\la_i \in \R$, $i=0,\cdots ,\ell_1$, are the mutually distinct 
real roots
and besides these there also the pairs of complex conjugate roots
$-\ka_m, -\bar{\ka}_m \in \C \setminus \R$, with the $\ka_m$ mutually 
distinct for $m=0,\cdots,\ell_2$. 
So from Theorem
\ref{polyid} we have
\begin{equation}\label{cxsum}
 1=\sum_{i=0}^{\ell_1} \Pr^N_{\lambda_i}[x] +
\sum_{m=0}^{\ell_2}(\Pr^N_{\kappa_m}[x] + \Pr^N_{\bar{\ka}_m}[x])~,
\end{equation} where we have made an obvious adaption of the notation. By
inspecting the formula there (i.e. \nn{genpolyid}) we see that
$\Pr^N_{\lambda_i}[x]$ is real, for $i=0,\cdots ,\ell_1$, and so is
each sum $(\Pr^N_{\kappa_m}[x] + \Pr^N_{\bar{\ka}_m})$, $m=0,\cdots
\ell_2$. We note also that in each $(\Pr^N_{\kappa_m}[x] +
\Pr^N_{\bar{\ka}_m})$ there is a common factor $P^{m\bar
m}:=P[x]/(x+\ka)^{p_{\ka_m}}(x+\bar{\ka})^{p_{\ka_m}}$.  Thus, in
summary, by combining conjugate factors we obtain a real identity in
$\mathbb{R}[x]$ of the form
\begin{equation}\label{cxsum2}
  1=\big(\sum_{i=0}^{i=\ell_1}Q_iP^i\big)+ 
  \big(\sum_{m=0}^{m=\ell_2}(Q_{m\bar m} P^{m\bar m})\big)~,
\end{equation} 
where the $Q_i$ and the $Q_{m\bar{m}}$ are obtained explicitly from
\nn{cxsum}, and where each term in the sum has polynomial degree less
than the degree of $P[x]$.  Thus from the general results of Section
\ref{gendecom} we obtain the following.

\begin{corollary}[A real version of Theorem \ref{inh}]
Let $\cV$ be a real vector space and $P = P[\cD]: \cV \to \cV$ an operator
polynomial in $\cD: \cV \to \cV$. Assume the complexification of $P$ 
factors as in \nn{Popformreal} where $\la_i \in \R$ and 
$\ka_j = \mu_j+i\nu_j \in \C \setminus \R$ with
$\mu_j,\nu_j \in \R$,
$\nu_j \not= 0$ and the $\la_i$'s and $\ka_j$'s are mutually distinct.
Then the null space $\cV_P$, for $P$, admits a canonical
and unique direct sum decomposition
\begin{equation} \label{sum:genR}
\cV_P=\bigoplus_{i=0}^{\ell_1} \cV_{\lambda_i} \oplus
      \bigoplus_{k=0}^{\ell_2} \cV_{\mu_k,\nu_k}
\end{equation}
where for each $k$ in the sum, $\cV_{\mu_k,\nu_k}$ is the solution
space for $(\cD^2 + 2\mu_j\cD + \mu_j^2+\nu_j^2)^{q_k}$.

Fixing $f \in \cV$, there is a 1-1 
relationship between solutions $u \in \cV$ of $Pu=f$ and solutions
$(u_0,\ldots,u_{\ell_1},u'_0,\ldots,u'_{\ell_2})$ of the problem
\begin{gather*}
  (\cD+\la_0)^{p_0}u_0=f, \ldots, (\cD+\la_{\ell_1})^{p_{\ell_1}}
  u_{\ell_1} = f, \\ (\cD^2 + 2\mu_0\cD + \mu_0^2 + \nu_0^2)^{q_0}
  u'_0 = f, \ldots, (\cD^2 + 2\mu_{\ell_2}\cD + \mu_{\ell_2}^2 +
  \nu_{\ell_2}^2)^{q_{\ell_2}} u'_{\ell_2} = f. 
\end{gather*}The mappings relating
  these are given by $F$ and $B$ in Theorem \ref{inhg} using
  \nn{cxsum2}.
\end{corollary}

\subsection{Operators polynomial in commuting endomorphisms}
\label{multivar}

We  now move to the general situation for this section. As above
let us write $\cV$ to denote a vector space over some field $\bF$.
Suppose that $\cD_i:\cV\to \cV$, $i=1,\cdots ,k$, are non-trivial
linear endomorphisms that are mutually commuting:
$\cD_i\cD_j=\cD_j\cD_i$ for $i,j\in\{1,\cdots ,k\}$ . We obtain a
commutative algebra $\bF[\bD]$ of consisting of those endomorphisms
$\cV\to \cV$ which may be given by expressions polynomial (with
coefficients in $\bF$) in the $\cD_i$.  We write $\bx =
(x_1,\ldots,x_k)$ for the multivariable indeterminate, and $\F[\bx]$
for the algebra of polynomials in the variables $x_1,\ldots,x_k$ over
the field $\F$.  Generalising the case of single variable polynomials,
there is a unital algebra epimorphism from $\bF[\bx]$ onto $\bF[\bD]$
given by formally replacing each variable $x_i$, in a polynomial, with
$\cD_i$.  

Given polynomials $P_0[\bx], P_1[\bx], \cdots ,P_{\ell}[\bx] \in \F[\bx]$
consider the product polynomial
\begin{equation} \label{Pmulti}
  P[\bx] = P_0[\bx] P_1[\bx] \cdots P_{\ell}[\bx].
\end{equation}
With $L=\{0,1,\cdots ,L\}$, we carry over, in an obvious way, the
labelling from Sections \ref{complexes} and \ref{aldec} via elements
of the power set $2^L$; products of the polynomial $P_i[\bx]$ are
labelled by the corresponding subset of $L$. For example for
$J\subseteq L$, $P_J[\bx]$ means $\prod_{j\in J}P_j[\bx]$, while
$P^J[\bx]$ mean $P_{L\setminus J}[\bx]$.

With a view to linking to the constructions above, we seek polynomials
$Q_J[\bx] \in \F[\bx]$, $J\in \al \subseteq 2^L$ satisfying the identity
\begin{equation} \label{alProjmulti}
  1 = \sum_{J\in \al} Q_J[\bx] P^J[\bx]~,
\end{equation}
or equivalently
\begin{equation} \label{almultieuclid}
  1 \in \langle P^J[\bx]~:~ J\in \al \rangle
\end{equation}
where $\langle .. \rangle$ denotes the ideal in $\F[\bx]$ generated by
the enclosed polynomials.
Via the polynomial analogue of Proposition \ref{Peucgen} 
we may
equivalently study the ``dual'' problem of finding sets $\beta \subseteq 2^L$
so that for each $I\in \beta$ we have 
\begin{equation}\label{duale}
1 \in \langle P_i[\bx]~:~ i\in I  \rangle~.
\end{equation}

We may use algebraic geometry to shed light on this problem.  Let us
 write $\cN(S[\bx])$ for the algebraic variety determined by the
 polynomial $S[\bx] \in \F[\bx]$ (i.e. $ \{ \bx \in \F^k \mid S[\bx]=0
 \}$) and put $\cN_I := \cN(P_I)$ and $\cN^J := \cN(P^J)$ for $I,J\in
 2^L$.  Clearly the condition \nn{duale} requires $\cap_{i\in I}\cN_i
 = \emptyset$, because $\cap_{i\in I}\cN_i \subseteq \cN(1) =
 \emptyset$.  Comparing to \nn{duale}, the condition $\cap_{i\in
 I}\cN_i = \emptyset$ is easier to verify, at least in simple cases,
 but it is generally weaker.  However this depends on the field. In
 particular, it follows from the (weak form of) Hilbert's
 Nullstellensatz (see e.g.\ \cite[Chapter 4, Theorem 1]{CLS}) that if
 $\F$ is algebraically closed then
$$  
1 \in \langle P_i[\bx]~:~ i\in I  \rangle
\Longleftrightarrow \cap_{i\in I}\cN_i
 = \emptyset ; 
$$ 
i.e.\ for $\F$ algebraically closed the condition 
$ 1 \in \langle P_i[\bx]~:~ i\in I  \rangle$ is equivalent to the polynomials
 $ P_i[\bx]~:~ i\in I $ having no common zero. (Note the previous
 display does not hold for $\F = \R$, e.g.\ take $P_i[x,y] = x^2+1$
 and $P_j[x,y] = y^2$.) Thus with the notation introduced at the start
 of this section we have the following.
\begin{theorem}\label{geomst}
For $P[\bD]\in \bF[\bD]$, ($\bD=(\cD_1,\cdots ,\cD_k)$) with $\bF$
algebraically closed,
\begin{equation}\label{polyD}
P[\bD]= P_0[\bD] P_1[\bD]\cdots  P_\ell[\bD]
\end{equation}
 is an algebraic dual $\beta$-decomposition of $P[\bD]$ if and only
if, for the polynomials $P_i[\bx]$ corresponding to the factors
$P_i[\bD]$, we have
$$ \cap_{i\in I}\cN_i
 = \emptyset, \quad \mbox{ for all } I\in\beta~.
$$

Generically for $\beta \in 2^L$,  such that for all $I\in \beta $, $|I|\geq k$,
 \nn{polyD} is a dual $\beta$-decomposition.  
\end{theorem}

Here, in an obvious way, we are using the term {\em algebraic dual
$\beta$-decomposition} to mean a dual $\beta$-decomposition that
arises from the analogous polynomial identities as discussed.  The
last statement holds because generically $ \cap_{i\in I}\cN_i$ has
codimension $|I|$.

The Theorem indicates immediately why one expects very strong results
in the case of operators polynomial in a single operator. In one
dimension algebraic varieties are generically disjoint. According to
the theorem the situation is not much weaker for operators polynomial
in several commuting operators. If we fix $k$ then still we may say
that generically operators polynomial in the operators $\cD_0,\cdots
\cD_k$ admit (algebraic) $\alpha$-decompositions.
It is clear that Theorem \ref{geomst} may be used to easily construct
examples of all varieties of $\al$-decompositions.

We note that the results here are perhaps suggested by the
 general ideas of algebraic invertibility developed in \cite{Gromov}
 and references therein. However explicit links with the development
 in that source  are currently far from clear.

\section{Symmetries}\label{sym}

\def\cW{\mathcal{W}} \def\cS{\mathcal{S}} \def\cE{\mathcal{E}}

Suppose that $P$ is a linear endomorphism of a vector space $\cV$,
over a field $\bF$. As above we write $\cV_P$ for the kernel of
$P$. Let us say that a linear map $S:\cV\to \cV$ is a {\em strong
symmetry} of $P$ if $S$ preserves each of the 
 eigenspaces of $P$.  For example, if a
$\cV$ endomorphism $S$ commutes with $P$, that is on $\cV$ we have
$[S,P]:=SP-PS=0$, then $S$ is a strong symmetry.  On the other hand
let us say that a linear operator $S:\cV_P\to \cV$ is a {\em weak
symmetry} of $P$ if $S$ has image in $\cV_P\subset \cV$. That is if
$S$ takes $P$-solutions to $P$-solutions. For example, if $S:\cV_P\to
\cV$ satisfies $PS=S'P$ for some linear operator $S':\cV\to \cV$, then
$S$ is a weak symmetry.  Evidently weak symmetries may be composed and
via this operation yield an algebra. Similarly for strong symmetries.

Given $P$ as above, let us write $\cW_P$ for the space of weak
symmetries of $P$.  In the case that $P$ admits an algebraic
decomposition $P=P_0P_1\cdots P_\ell$ (as in section \ref{gendecom})
then we obtain a corresponding decomposition of $\cW_P$, as a vector
space. First one further item of notation. Let us write $\cW_{ij}$ for
the vector space of linear homomorphisms $H:\cV_{j}\to \cV_{i}$ where,
recall, $\cV_i$ is the null space of $P_i$. Here we carry
over notation from  Section \ref{gendecom}.
\begin{theorem}\label{genweak}
For $P:\cV\to \cV$, with \nn{compP} giving an algebraic decomposition, we have
a canonical vector space decomposition,
$$
\cW_P\cong \oplus_{i,j=0}^{i,j=\ell}\cW_{ij}.
$$ 
\end{theorem}
\begin{proof}
For $H \in \cW_{ij}$ w obtain an element in $\cW_P$ by forming 
$H\circ Proj_j$. This is inverted by the map taking arbitrary 
 $S\in \cW_P$ to  the composition
$$
\Proj_i\circ S\circ \Proj_j
$$
 in $\cW_{ji}$. 
\end{proof}
Note that for $H_{jk}\in \cW_{jk}$ and $H_{ij}\in \cW_{ij}$ we
have $H_{ij}\circ H_{jk}\in \cW_{ik}$. Thus, identifying $\cW_P$ with
$\oplus \cW_{ij}$ via the isomorphism in the Theorem, we see that for
each $i=0,1,\cdots ,\ell$, $\cW_{ii}$ is a subalgebra of
$\cW_P$. Evidently the algebra structure of $\cW_P$ arises from that
of these subalgebras plus the interlacing introduced by the spaces of
homomorphisms $\cW_{ij}$, where $i$ and $j$ are distinct. 
Overall,
understanding the algebraic structure of $\cW_P$ is reduced to
understanding the spaces $\cW_{ij}$.

Now suppose that $\bF$ is an algebraically closed field and $P$ is any
polynomial in $\cD$. Recall from Corollary \ref{spectralthm} that for
a given $\mu\in \bF$ the corresponding $P$-eigenspace (for simplicity
of discussion we will allow this to be possibly trivial) $\cV^\mu$
decomposes into a direct sum
$\cV^\mu=\oplus_{i=1}^{k_\mu}\cV_{\lambda_i}$ where the
$\cV_{\lambda_i}$ are generalised eigenspaces for $\cD$. 
Evidently we have the following observation.
\begin{proposition}\label{sttry}
If $P:\cV\to \cV$ is a linear operator non-trivially polynomial in $\cD$
and $S:\cV\to \cV$ preserves all generalised eigenspaces for $\cD$,
then $S$ is a strong symmetry for $P$.
\end{proposition}
\noindent So for example any polynomial in $\cD$ (viewed as a linear operator
$\cV\to \cV$) is a strong symmetry. 

The conditions in the Proposition are obviously too strict to generate
all strong symmetries in general. 
It would be
interesting to understand the precise relationship between strong
symmetries for operators $P$, as in the Proposition, and the eigenspace
information for $\cD$. As a passing note we make a final observation
in this direction.  It is clear that if we fix $\mu$ in $\bF$ then the
restriction to $\cV^\mu$ of the linear maps $S:\cV\to \cV$ that
preserve $\cV^\mu$ yields a space $\tilde{\cW}_{P-\mu}$ which is
defined in the same way as the space of weak symmetries for the
operator $P-\mu$, except that it consists of maps $\cV\to \cV$ (the
domain is not taken to be the solution space). 
Thus this may be analysed as was done for $\cW_P$ above.  The
situation is rather simple in lower degree cases. For example, the
following proposition describes strong symmetries explicitly for $P$
of degree $2$.
\begin{proposition}
  Let $P = (\cD+\la_1)(\cD+\la_2)$, with $\la_1,\la_2 \in \mathbb{F}$
  (not necessarily distinct) and where $\cD: \cV \to \cV$ is a linear
  operator.  For $\xi\in \bF$, denote by $\cV_{\xi}^k$ the solution
  space of $(\cD+\xi)^k$ for $k \in \N$. Then $S: \cV \to \cV$ is a
  strong symmetry of $P$ if and only if the following three conditions
  hold: \\
  (i) if $-\xi_0 := -\frac{\la_1+\la_2}{2} \in \mathrm{Spec}\, D$
  then $S$ preserves $\cV_{\xi_0}^2$ \\
  (ii) if $-\xi, -(\la_1+\la_2-\xi) \in \mathrm{Spec}\, D$, $\xi \not=
  \xi_0$ then $S$ preserves
  $\cV_{\xi}^1 \oplus \cV_{\la_1+\la_2-\xi}^1$ \\
  (iii) if $-\xi \in \mathrm{Spec}\, D \wedge -(\la_1+\la_2-\xi)
  \not\in \mathrm{Spec}\, D$ then $S$ preserves $\cV_{\xi}^1$.
\end{proposition}

\begin{proof}
Consider the decomposition of $P-\mu$ to irredicibles, i.e.\ 
$$ P-\mu = (\cD+\la_1)(\cD+\la_2) - \mu = (\cD+\xi_1)(\cD+\xi_2) $$
where $\mu \in \mathbb{F}$ and $\xi_1,\xi_2 \in \mathbb{F}$
are not necessarily distinct. 
 Then clearly $\xi_1+\xi_2 = \la_1+\la_2$
and any pair $\xi_1,\xi_2$ such that $\xi_1+\xi_2 = \la_1+\la_2$
satisfies the previous display for some $\mu \in \mathbb{F}$.
Thus the strong symmetries are precisely linear mappings preserving
the solution space of $(\cD+\xi_1)(\cD+\xi_2)$ for every $\xi_1,\xi_2$
such that $\xi_1+\xi_2 = \la_1+\la_2$. Using Theorem \ref{pf0}, the 
proposition follows.
\end{proof}

\section{Conformal Laplacian operators and Einstein manifolds} \label{eins}

On a smooth Riemannian or pseudo-Riemannian manifold $(M^n,g)$ let
us write $d$ for the exterior derivative and $\delta$ for its formal
adjoint (as in e.g.\ \cite{Besse}). On the space of smooth
$k$-forms $\Lambda^k$ the form Laplacian is given by 
$\Delta = \delta d+ d\delta $. Consider the operator
$$
Z=\Delta -\lambda^2
$$
where $0 \neq \lambda\in {\Bbb C}$. We may extend this to an operator on 
$\Lambda^*$, the space of all smooth differential forms. 
Thus we have $\Delta =
D^2$ where $D$ is the Dirac operator $ d+\delta$,  hence 
$$
Z=(D+\lambda)(D-\lambda)
$$ 
on $\Lambda^*$.
Thus fixing $f\in \Lambda^*$, solutions $u\in \Lambda^0$ 
of the problem $Z u=f$  are in 1-1 correspondence with solutions 
$(u_+,u_-)\in (\Lambda^*)^2$ of the problem 
$$
(D+\la)u_+=f \quad (D-\la)u_-=f~,
$$
where we view $\Lambda^0\subset \Lambda^*$.
The map from $u$, solving $Z u=f$, to a solution
of the display is 
\begin{equation}\label{Fs}
u\mapsto \big( (D-\la)u,(D+\la)u \big)~,
\end{equation}
while the inverse is 
\begin{equation}\label{Bs}
(u_+,u_-)\mapsto \frac{1}{2\la}(u_--u_+).
\end{equation}

Using the grading of forms by degree, we may apply these tools to $Z$
as an operator on functions $\Lambda^0$. It is easily seen that
\nn{Bs} specialises to a map from $\Lambda^0\to(\Lambda^0\oplus
\Lambda^1)^2 $, inverted by \nn{Bs} as a map $(\Lambda^0\oplus
\Lambda^1)^2\to \Lambda^0$.  Fixing $f\in \Lambda^0$, this gives a 1-1
relationship between functions $u$ solving $Zu=f$ and pairs
$(u_+,u_-)\in (\Lambda^0\oplus \Lambda^1)^2$ solving $ (D+\la)u_+=f$
and $(D-\la)u_-=f$. In fact, once again using the grading of forms,
one sees that the second order equation $Z u=f $ is in fact equivalent
to either one of the first order equations $ (D+\la)u_+=f$
or $(D-\la)u_-=f$.

Operators of the form of $Z$ arise naturally in Riemannian geometry.
The conformal Laplacian 
$Y:\Lambda^0\to \Lambda^0$  is given by the formula
$$
Y=\Delta+\frac{n(n-2)}{4n(n-1)}{\rm Sc}
$$ 
where ${\rm Sc}$ is the scalar curvature.  Thus this is of the same
form as $Z$ on manifolds where ${\rm Sc}$ is constant and non-zero.  

A 
class of constant scalar curvature  manifolds are the Einstein structures.  A
pseudo-Riemannian manifold $(M,g)$ is said to be Einstein if its Ricci
curvature is proportional to the metric (all structures will be taken
to be smooth).  We refer the reader to \cite{Besse} for background on
the meaning of these statements and the importance of Einstein
structures.  The conventions below follow \cite{GoEinst} except that we
will use the ``positive energy'' Laplacian $\De$ as above (it may be 
also given as $\Delta =\nabla^*\nabla$,
where $\nabla$ is the Levi-Civita connection and $\nabla^*$ its formal
adjoint.  We assume the dimension of $M$ to be at least 3.
The GJMS conformal Laplacians of \cite{GJMS} are in general given by
extremely complicated formulae, see \cite{GoPetLap}. However on
conformally Einstein manifolds we may choose an Einstein metric $g$.
Then the order $2k$ GJMS operator may be viewed as an operator
$P_{k}:\Lambda^0\to\Lambda^0$ and the formulae for these may be
simplified dramatically.  On Einstein $n$-manifolds the $P_k$ is given
by \cite{GoEinst}
\begin{equation}\label{prodformula}
P_{k} = \prod_{i=1}^k(\De +c_i {\rm Sc}),
\end{equation}
where $c_i= (n+2i-2)(n-2i)/(4n(n-1))$ and ${\rm Sc}$ is the scalar
curvature, that is the metric trace of the Ricci curvature.  (For the
standard sphere as a special case the formula \nn{prodformula} was
known to Branson \cite{BransonSoeul}.)  On even manifolds the GJMS
operators exist only up to order $n$. However for conformally Einstein
structures it is shown in \cite{GoEinst} that, in a suitable sense,
the family extends to all even orders. So for our current purposes for
any $k\in \mathbb{Z}_{>0}$ we term the operator \nn{prodformula} a
GJMS operator. (We should also note that in line with our conventions
for the sign of the Laplacian, the GJMS operator $P_k$ as above is
$(-1)^k$ times the corresponding operator in \cite{GoEinst}).

Since the scalar curvature ${\rm Sc}$ is necessarily constant on
Einstein manifolds it follows that $P_k$ is polynomial in $\De$ and so
we may immediately apply the results above to relate the null space of
$P_k$ with the generalised eigenvalues of the Laplacian.  In the
setting of compact manifolds of Riemannian signature it was noted in
\cite{GoEinst} that we have such information via standard Hodge theory
(or one could use functional calculus). The gain here is that we
obtain related information in any signature and without any assumption
of compactness.

The left (i.e.\ $i=1$) factor in the expansion \nn{prodformula} is in
fact the conformal Laplacian $Y$ which plays a central role in
spectral theory. So let us instead rephrase the Theorem 1.2 from 
\cite{GoEinst} in terms of this. 
\begin{theorem}
On a pseudo-Riemannian $n$-manifold with Einstein metric, 
 the order $2k$ GJMS operator is given by
\begin{equation}\label{prodformulaII}
P_k = \prod_{i=1}^k(Y + b_i {\rm Sc}),
\end{equation}
where $b_i= \frac{i(1-i)}{n(n-1)}$.
\end{theorem}

Note that when ${\rm Sc}\neq 0$ the 
scalars $b_i {\rm Sc}$
are mutually distinct.  Thus from Theorem
\ref{pf0:gen}, and writing $\cN(P_k)$ for the null space of $P_k$ as
an operator on smooth real valued functions, we have the following.
\begin{theorem}\label{maina}
On a pseudo-Riemannian Einstein $n$-manifold with ${\rm Sc}\neq 0$ the
null space of $P_k$ has a direct sum decomposition
$$
\cN(P_k)= \oplus_{i=1}^k \cN_i(Y)~,
$$
where $\cN_i(Y)$ is the eigenspace for $Y$ with eigenvalue $-b_i {\rm Sc}$.
\end{theorem}
\noindent Of course the machinery implies in the case of ${\rm Sc}=
0$, but in this case the result is obvious: the null space is a
generalised eigenspace for $Y$ with generalised eigenvalue 0, that is
$\cN(P_k)= \cN(Y^k)$.  In all cases the projection $\cN(P_k) \to\cN_i(Y)
$ is given by \nn{projform:gen}.  Similarly, the eigenspectrum of $P_k$ is
determined by Corollary \ref{spectralthm}.

\begin{theorem}\label{specPk}
On a pseudo-Riemannian Einstein $n$-manifold,  $(\mu, f)$ is an
eigenvalue, eigenfunction pair for the GJMS operator $P_k$ if and only
if for some $m\in\{1,\cdots ,k\}$
$$
f=f_1+\cdots+f_m~, \quad \quad 0\neq f_i,\quad i=1,\cdots , m,
$$
where, for each $i\in\{1,\cdots , m\}$, 
$(Y-\la_i)^{p_i} f_i =0 $ and $\la_i$ is a multiplicity $p_i$ solution of 
of the polynomial equation $(P_k-\mu)[x]=0$. (Here we consider $P_k$ as the 
polynomial in $Y$, i.e.\ given by \nn{prodformulaII}.)
\end{theorem}
The inhomogeneous problems yield the obvious simplification to second
order problems.
\begin{proposition}\label{inhPk}
On a pseudo-Riemannian Einstein $n$-manifold, 
the inhomogeneous problem $P_k u=f$, for the GJMS operator $P_k$,
is equivalent to the second order problem 
$$
(Y+b_1 {\rm Sc})u_1=f, \cdots , (Y+b_k {\rm Sc})u_k=f~.
$$
From a solution $(u_1,\cdots ,u_k)$ of this problem we obtain,
using $b_i= \frac{i(i-1)}{n(n-1)}$, the solution $u$ of $P_k u=f$ as
$$ u = \bigl( \frac{n(n-1)}{{\rm Sc}} \bigr)^{k-1}
   \sum_{i=1}^{k} \bigl[ \prod_{i\neq j =1}^{j=k} 
   \frac{1}{(j-i)(j+i-1)} \bigr] u_i. $$
\end{proposition}
\noindent In fact in odd dimensions and also in even dimensions $n$
for the operators $P_{k\leq n/2}$ we may further reduce to first order
problems using the ideas at the start of this section. Via different
Dirac operators there are variations on this outcome.

\subsection{Differential Weak symmetries}\label{natweak}

 It is clear that in any special setting the general idea of
symmetries may be tuned somewhat. In particular, we shall do this for
differential operators on pseudo-Riemannian manifolds. Suppose that
now $\cV$ is a space of smooth sections of some vector bundle over a
pseudo-Riemannian manifold and $P:\cV\to \cV$ is a differential
operator.  Then we shall say that a weak symmetry $S$ of the
differential operator $P$ is differential if is given by a
differential operator on $\cV$. That is $S$ is differential weak
symmetry of $P$ means that it is a differential operator $S:\cV\to\cV$
such that it preserves the solution space of $P$.  (This is slightly
different from Section \ref{sym} where we defined weak symmetries only
on the solution space of $P$.)  Since the composition of differential
operators yields a differential operator the differential weak
symmetries form a subalgebra of the weak symmetries for $P$.  Similar
ideas apply to strong symmetries which may also be required to be
differential.  The key point is that provided the
projection operators $\Proj_i$ (from Corollary \ref{maincor}) are
differential then the general results from section \ref{sym} carry over
functorially to this category. 

In particular we illustrate this in the setting as above. Here we take
$\cV$ to be the space of smooth functions $\cE$ on an Einstein
manifold $M$ (of dimension at least 3).  Let us write $\cW^{P_k}_{ij}$
for the space of linear differential operators $S:\cE\to \cE$ with the
property that, upon restriction to $\cN_j$, $S$ takes values in
$\cN_i(Y)$, that is $S:\cN_j\to \cN_i(Y)$.  The differential operators
in $S$ map between eigenspaces of the conformal Laplacian $Y$. From
Theorem \ref{pf0:gen} and Theorem \ref{genweak} we deduce the
following.
\begin{theorem}\label{GJMSweak}
On a pseudo-Riemannian Einstein $n$-manifold with ${\rm Sc}\neq 0$,
the  space $\cW_{P_k}$ 
of differential weak symmetries of the order $2k$ GJMS operator $P_k$
has a canonical vector space decomposition,
$$
\cW_{P_k}\cong \oplus_{i,j=0}^{i,j=k}\cW^{P_k}_{ij}.
$$ 
\end{theorem}

An obvious specialisation is to consider conformally flat spaces and
locally (i.e.\ on a contractible manifold). Since the GJMS operators
are conformally invariant, their solution spaces are conformally
stable and one may study these by choosing a conformal scale that is
congenial for the problem. For a current purposes a scale that
achieves a constant non-zero curvature is ideal since then (on such
Einstein structures) Theorem \ref{maina} applies. In particular we may
apply Theorem \ref{GJMSweak} to study this conformal problem.  In the
setting of Euclidean space, Eastwood and Eastwood-Leistner
\cite{MikEsym,EL} have studied the ``higher symmetries'' of the
Laplacian and its square. These are differential weak symmetries $S$
with the property that (in a choice of conformal scale) $PS=S'P$ where
$S':\cE\to \cE$ is a differential operator.  In this flat setting the
Laplacian agrees with the Yamabe operator while the square of the
Laplacian is the order 4 GJMS operator (which is usually termed the
Paneitz operator).  Since their theory is essentially conformal it
should be an interesting direction to carry their results for the
square of the Laplacian, in \cite{EL}, onto a constant curvature
conformally flat space and then relate these to our observations
above. Our tools above provide an alternative approach to such higher
order problems and also provide a route for studying the related
questions on general conformally Einstein manifolds.


\begin{thebibliography}{XX}

\bibitem{Besse} A. L. Besse, ``Einstein manifolds'', Springer-Verlag, 
Berlin, 1987. xii+510

\bibitem{BransonSoeul} T. Branson, ``The Functional Determinant'',
Global Analysis Research Center Lecture Note Series, Number 4,
Seoul National University (1993).

\bibitem{tomsharp} T. Branson, \textit{Sharp
  inequalities, the functional determinant, and the complementary
  series}.  Trans.\ Amer.\ Math.\ Soc.\ {\bf 347} (1995) 3671--3742.

\bibitem{BrGodeRham} T.\ Branson, and A.\ R.\ Gover, \idx{Conformally
  invariant operators, differential forms, cohomology and a
  generalisation of $Q$ curvature}, Comm.\ Partial Differential
  Equations, {\bf 30} (2005), 1611 - 1669.

\bibitem{CLS} D. Cox, J. Little, D. O'Shea, ``Ideals, varieties,
and algorithms. An introduction to computational algebraic geometry
and commutative algebra.'' Second edition. Undergraduate Texts in
Mathematics. Springer-Verlag, New York, 1997. xiv+536 pp.

\bibitem{dirac} P.A.M. Dirac,  {\em Wave equations in conformal
space}.  { Ann.\ of Math.\/} {\bf 37}, (1936) 429--442.

\bibitem{DM} Z.\ Djadli and A.\ Malchiodi, {\em Existence of conformal
    metrics with constant $Q$-curvature}.
Preprint math.AP/0410141,  http://www.arxiv.org

\bibitem{MikEsym} Michael Eastwood, {\em Higher symmetries of the
Laplacian}, Ann.\ of Math.\ {\bf 161} (2005), 1645--1665.

\bibitem{EL} Michael Eastwood, and Thomas Leistner, {\em Higher
Symmetries of the Square of the Laplacian}, preprint math.DG/0610610.

\bibitem{Esc} J.\ Eschmeier, {\em Local properties of Taylor's
    analytic functional calculus}, Invent.\ Math.\ {\bf 68} (1982),
  103--116.

\bibitem{FeffGr01} C.\ Fefferman, C.R.\ Graham,  {\em $Q$-curvature
    and Poincar\'{e} metrics},  Math.\ Res.\ Lett.\ {\bf 9}, 139-151 (2002).

\bibitem{GoEinst} A.R.\ Gover, {\em Laplacian operators and
    Q-curvature on conformally Einstein manifolds}, Mathematische
  Annalen, {\bf 336} (2006), 311--334.


\bibitem{GoGr} A.R. Gover, C.R. Graham, \textit{CR Invariant
    Powers of the sub--Laplacian} J. Reine Angew. Math. \textbf{583}
  (2005), 1--27.  

\bibitem{GoPetLap} A.R.\ Gover and L.J.\ Peterson, \textit{Conformally
invariant powers of the Laplacian, Q-curvature, and tractor calculus}.
Commun.\ Math.\ Phys.\ {\bf 235} (2003) 339--378.


\bibitem{GJMS} C.R.\ Graham, R.\ Jenne, L.J.\ Mason, G.A.\ Sparling,
  \textit{Conformally invariant powers of the Laplacian, I:
    Existence}. J. London Math. Soc.\ \textbf{46}, (1992) 557--565.

\bibitem{GrZ2} C.R.\ Graham, M.\ Zworski, {\em Scattering matrix in
    conformal geometry},  Invent.\ Math., {\bf 152}  (2003), 89--118.

\bibitem{Gromov} M.\ Gromov, Partial differential relations,
Ergebnisse der Mathematik und ihrer Grenzgebiete (3), 9
Springer-Verlag, Berlin, 1986. x+363 pp.


\bibitem{MCrelle} V.\ M\"uller, {\em Local behaviour of the polynomial
  calculus of operators},  J.\ Reine Angew.\ Math.\  {\bf 430} (1992), 61--68.

\bibitem{M} V.\ M\"uller, Spectral theory of
  linear operators and spectral systems in Banach algebras. Operator
  Theory: Advances and Applications, 139. Birkhäuser Verlag, Basel,
  2003, x+381 pp.


\bibitem{Pan} S.\ Paneitz, {\em A quartic conformally covariant
differential operator for arbitrary pseudo-Riemannian manifolds}.
Preprint (1983).

\bibitem{schoen} R.\ Schoen, \idx{Conformal deformation of 
a Riemannian metric to constant scalar curvature},
J. Differential Geom. 20 (1984), no. 2, 479--495

\bibitem{JoSth} J.\ \v Silhan, Invariant operators in conformal
geometry, PhD thesis, University of Auckland, 2006.

\bibitem{T} J.L.\ Taylor, \idx{A joint sepctrum for several commuting
operators}, J.\ Funct.\ Anal.\ {\bf 6}, (1970), 172--191.

\end{thebibliography}
\end{document}